\documentclass{amsart}
\usepackage{amssymb,euscript,amsmath, mathrsfs}
\usepackage[dvips]{graphicx}
\usepackage[dvips]{color}

\newcounter{ENUM}
\newcommand{\itm}{\item}
\newenvironment{ilist}{\renewcommand{\theENUM}{\roman{ENUM}}\renewcommand{\itm}{\addtocounter{ENUM}{1}\item[(\theENUM)]}\begin{itemize}\setcounter{ENUM}{0}}{\end{itemize}}
\newenvironment{Ilist}{\renewcommand{\theENUM}{\Roman{ENUM}}\renewcommand{\itm}{\addtocounter{ENUM}{1}\item[(\theENUM)]}\begin{itemize}\setcounter{ENUM}{0}}{\end{itemize}}

\newcommand{\margh}[1]{}

\input xy
\xyoption{all}
\CompileMatrices

\def\A{{\mathbb A}}
\def\P{{\mathbb P}}

\def\G{{\mathbb G}}

\def\F{{\mathscr F}}

\def\L{{\mathscr L}}

\def\E{{\mathscr E}}
\def\O{{\mathscr O}}

\def\sP{{\mathscr P}}

\def\sN{{\mathscr N}}
\def\cG{{\mathcal G}}

\def\cLG{\mathcal{LG}}

\def\g{{\mathfrak g}}

\def\m{{\mathfrak m}}

\def\Tor{\operatorname{Tor}}

\def\sep{\operatorname{sep}}

\def\Gal{\operatorname{Gal}}

\def\Spec{\operatorname{Spec}}
\def\Pic{\operatorname{Pic}}

\def\univ{\operatorname{univ}}

\def\rk{\operatorname{rk}}
\def\im{\operatorname{im}}

\def\refn{\operatorname{ref}}
\def\red{\operatorname{red}}

\def\LG{\operatorname{LG}}

\def\EH{\operatorname{EH}}

\numberwithin{equation}{section}
\newtheorem{thm}{Theorem}[section]
\newtheorem{prop}[thm]{Proposition}
\newtheorem{lem}[thm]{Lemma}
\newtheorem{cor}[thm]{Corollary}

\theoremstyle{definition}
\newtheorem{defn}[thm]{Definition}
\newtheorem{ques}[thm]{Question}
\newtheorem{ex}[thm]{Example}
\newtheorem{sit}[thm]{Situation}
\theoremstyle{remark}
\newtheorem{notn}[thm]{Notation}
\newtheorem{rem}[thm]{Remark}
\newtheorem{warn}[thm]{Warning}

\begin{document}
\title{A Limit Linear Series Moduli Scheme}
\author{Brian Osserman}
\begin{abstract}
We develop a new, more functorial construction for the basic theory
of limit linear series, which provides a compactification of the
Eisenbud-Harris theory, and shows promise for generalization to 
higher-dimensional varieties and higher-rank vector bundles.
We also give a result on lifting linear series from characteristic $p$ 
to characteristic $0$. In an appendix, in order to obtain the necessary
dimensional lower bounds on our limit linear series scheme we develop a
theory of ``linked Grassmannians;'' these are schemes parametrizing
sub-bundles of a sequence of vector bundles which map into one another under
fixed maps of the ambient bundles.
\end{abstract}
\thanks{This paper was partially supported by a fellowship from the 
Japan Society for the Promotion of Science.}
\maketitle

\section{Introduction}

The Eisenbud-Harris theory of limit linear series of \cite{e-h1} is a
powerful tool for degeneration arguments on curves, with applications to the
Kodaira dimension of moduli spaces of curves, and analysis of Weierstrass
points on curves, as well as new arguments for results such as
the Gieseker-Petri theorem. In this paper, we give a new construction for 
limit linear series, very much in the spirit of Eisenbud and Harris' theory, 
but more functorial in nature, and involving a substantially new approach 
which appears better suited to generalization to higher-dimensional
varieties and higher-rank vector bundles. The application of the theory of
limit linear series in positive characteristic is fundamental to \cite{os7};
we should remark that we do not see any obstructions to the original 
construction of Eisenbud and Harris working in
characteristic $p$, but the independence of characteristic is more
transparent in the functorial setting. We begin with an overview of the
basic ideas of limit linear series; for those unfamiliar with linear series,
the actual definitions and notation are all recalled below.

While our main theorem (see Theorem \ref{grd-main}) is too technical to 
state in an introduction, we can 
outline the main concepts involved. The basic idea of limit linear series is 
to analyze how linear series behave as a family of smooth curves $X/B$ 
degenerates to a nodal curve $X_0$; a key distinguishing feature of the
theory is that rather than standard deformation-theoretic techniques to
obtain results from the degeneration, a simple dimension count on the
special fiber produces results immediately. 

More specifically, recall 
that a proper, geometrically reduced and connected
nodal curve with smooth components is said to be of {\bf compact type} if
the dual graph is a tree, or equivalently if the (connected component of
the) Picard scheme is proper. 
Now, if $X_0$ is not of compact type, line bundles on the smooth curves 
may not limit to a line bundle on the nodal fiber, as the Picard scheme of 
the family (and specifically of the nodal fiber) will not be proper. On the 
other hand, if the nodal fiber is reducible, limiting line bundles will 
exist, but will
not be unique, as one can always twist by one of the components of the
reducible fiber to get a new line bundle, isomorphic away from the nodal
fiber to the original one. However, this turns out to be the only ambiguity.
To explain the approach to this issue, we consider for simplicity the case
of $\g^r_d$'s where the family $X/B$ has smooth general fiber, and $X_0$ 
consists of two smooth components $Y$ and $Z$, meeting at a single node $P$. 
Given a line bundle of degree $d$ on $X$, we will say it has degree
$(i,d-i)$ on $X_0$ if it restricts to a line bundle of degree $i$ on $Y$ and 
degree $d-i$ on $Z$.
Eisenbud and Harris approached the problem by considering the linear
series obtained by looking at the two possible limit line bundles obtained
by requiring degrees $(d,0)$ and $(0,d)$ on $X_0$. Since the degree $0$ 
components cannot contribute anything to the
space of global sections chosen for the $\g^r_d$, this is equivalent to
specifying a $\g^r_d$ on each of $Y$ and $Z$; they showed that
if such a pair arises as a limit of $\g^r_d$'s from the smooth fibers, it will 
satisfy the ramification condition
\begin{equation}a_i^Y(P)+a_{r-i}^Z(P) \geq d,\end{equation}
for all $i$, where $\{a_i^Y(P)\}_i$ and $\{a_i^Z(P)\}_i$ are the vanishing
sequences on $Y$ and $Z$ at the node. Eisenbud and Harris refer to such 
pairs on a nodal 
fiber as {\bf crude limit series}, and when the inequality is replaced by 
an equality, as {\bf refined limit series}.

Eisenbud and Harris' moduli scheme construction requires restriction to 
refined limit series, and as such is not generally proper, and is also
necessarily disconnected, being constructed as a disjoint union over the
different possible ramification indices at the nodes. Moreover, the
necessity to specify ramification indices makes it unsuitable for
generalizing from curves to higher-dimensional varieties. The basic idea of
our construction is to remember not just the line bundles of degree $(d,0)$
and $(0,d)$ on $X_0$, but also the $d-1$ line bundles of degree $(i,d-i)$ 
that lie in between. One can then replace the ramification condition with 
a simpler compatibility condition on the corresponding spaces
of global sections, yielding a very functorial approach to constructing the
moduli scheme. Further, one can show a high degree of compatibility with
Eisenbud and Harris' construction: in particular, for a curve (of compact
type) over a field, our construction contains the Eisenbud-Harris version as
an open subscheme. 

We begin in Section \ref{s-grd-background} with a review of the basics 
of linear series, but in arbitrary characteristic. In Section
\ref{s-grd-families} we give the precise conditions on the families of
curves we will consider, and show that such families may be contructed as
necessary. In Section \ref{s-grd-functor} we define the limit linear series
functors we will consider, and our main theorem, the representability of
these functors, is proved in Section \ref{s-grd-represent}; we conclude with
corollaries as in Eisenbud and Harris on smoothing linear series from the 
special fiber when the dimension is as expected, including in the cases of
positive and mixed characteristic. We compare our
theory to that of Eisenbud and Harris in Section \ref{s-grd-compare}, and
conclude with some further questions in Section \ref{s-grd-further}.
Finally, in Appendix \ref{s-grd-lg} we develop a theory of linked
Grassmannian schemes, which parametrize collections of sub-bundles of a
sequence of vector bundles linked together by maps between the bundles; this
is used in the construction of the limit linear series scheme in the main
theorem, and in particular to obtain the necessary lower bound on its
dimension.

The work here is of course entirely inspired by Eisenbud and Harris' original
construction in \cite{e-h1}. Attempts to generalize this theory have thus
far been sparse, but include, for instance, work of Esteves \cite{es1} and 
of Teixidor i Bigas \cite{te1} to generalize to certain curves not of 
compact type and higher-rank vector bundles respectively.

The contents of this paper form a portion of the author's 2004 PhD thesis at
MIT, under the direction of Johan de Jong.

\section*{Acknowledgements}

I would like to thank Johan de Jong for his tireless and invaluable
guidance. I would also like to thank Jason Starr, Joe Harris, Max
Lieblich, and Steve Kleiman for their many helpful discussions. 

\section{Linear Series in Arbitrary Characteristic}\label{s-grd-background}

Before getting into the technical definitions related to the central
construction, we begin with a few preliminary definitions and lemmas in the
case of a smooth, proper, geometrically integral curve $C$ of genus $g$, 
over a field $k$ of any characteristic.

First, recall:

\begin{defn}If $\L$ is a line bundle of degree $d$ on $C$, and $V$ an 
$(r+1)$-dimensional subspace of $H^0(C, \L)$, we call the pair $(\L, V)$ a 
a {\bf linear series} of degree $d$ and dimension $r$ on $C$, or a $\g^r_d$.
Given $(\L,V)$ a $\g^r_d$ on $C$, and a point $P$ of $C$, there is a unique
sequence of $r+1$ increasing integers $a^{(\L,V)}_i(P)$ called the {\bf
vanishing sequence} of $(\L,V)$ at $P$, given by the orders
of vanishing at $P$ of sections in $V$. We also define
$\alpha^{(\L,V)}_i(P) := a^{(\L,V)}_i(P)-i$, the {\bf ramification sequence}
of $(\L,V)$ at $P$. $(\L, V)$ is said to be {\bf unramified} at $P$ if all 
$\alpha^{(\L,V)}_i(P)$ are zero; otherwise, it is {\bf ramified} at $P$.
\end{defn}

\begin{warn}Since the ramification and vanishing sequences are equivalent
data, we tend to refer to conditions stated in terms of either one simply as
``ramification conditions.'' We will also drop the $(\L,V)$ superscript or
replace it as appropriate, particularly when we have a linear series on
each component of a reducible curve, when we will tend to simply use the
component to indicate which series we are referring to.
\end{warn}

The following definitions, being tailored to characteristic $p$, may be less 
standard:

\begin{defn}We say a linear series $(\L,V)$ on $C$ is {\bf separable} if it 
is not 
everywhere ramified. Otherwise, it is {\bf inseparable}. At a point $P$, 
we say that $(\L,V)$ is {\bf tamely ramified} if the characteristic is $0$
or if the vanishing orders $a_i(P)$ are maximally distributed mod 
$p$ (in particular, this holds at any unramified point).
Otherwise, we say that $(\L,V)$ is {\bf wildly ramified} at $P$. 
\end{defn}

The following result is a characteristic-$p$ version of a standard Pl\"ucker 
formula, whose proof simply adapts standard techniques:

\begin{prop}\label{grd-plucker}Let $C$ be a smooth, proper, geometrically 
integral curve of genus $g$ over a field $k$, and $(\L, V)$ a $\g^r_d$ on 
$C$. Then either $(\L, V)$ is inseparable, or we have the inequality 
$$\sum _{P \in C} \sum _i \alpha_i(P) \leq (r+1)d+\binom{r+1}{2}(2g-2).$$ 
Furthermore, this will be an equality if and only if $(\L,V)$ is everywhere 
tamely ramified; in particular, in this case inseparability is impossible.
\end{prop}

\begin{proof}We simply use the argument of \cite[Prop. 1.1]{e-h4}.
Even though it is intended for characteristic $0$, the proof follows through
equally well in characteristic $p$ for our modified statement, noting that
their ``Taylor expansion'' map is defined independent of characteristic, and
their formulas then hold on a formal level. Indeed, their
argument shows that if $(\L,V)$ induces a non-zero section 
$s(\L,V)$ of $\L^{\otimes r+1} \otimes 
(\Omega^1_C)^{\otimes \binom{r+1}{2}}$, we get the desired
inequality, with equality if and only if the determinant of their Lemma 1.2 
is non-zero at all $P$ (where, as in the proof of the proposition, $X_j :=
\alpha^{(\L,V)}_j(P)$). In fact, if this determinant is non-zero anywhere,
we see also that $s(\L,V)$ has finite order of vanishing at that point, and
cannot be the zero section. Next, their same lemma shows that their 
determinant will be non-zero at a point $P$ if and only if $(\L, V)$ is tamely 
ramified at $P$. This means that if we show that inseparability corresponds
precisely to having $s(\L,V)=0$, we are done. But
this also follows trivially, since on the one hand any unramified point is in 
particular tamely ramified, and will in fact give a non-vanishing point of
$s(\L,V)$, and on the other hand, if $s(\L,V)$ is non-zero, we have
seen that we can get only finitely many ramification points.
\end{proof}

Note that because vanishing sequences are bounded by $d$, if
$d<p$, then wild ramification is not possible, so the previous proposition
immediately implies:

\begin{cor}Wildly ramified or inseparable linear series of degree $d$ are 
only possible when $d \geq p$.
\end{cor}

Finally, we have the notation:

\begin{defn}\label{grd-expected} Given $n$ points $P_i$ and $n$ 
ramification sequences $\alpha^i = \{\alpha^i_j\}_j$, 
we write $\rho :=\rho(g,r,d; \alpha^i):=(r+1)(d-r)-rg-\sum_{i,j} \alpha^i_j$.
This is the {\bf expected dimension} of linear series of degree $d$ and 
dimension $r$ on a curve of genus $g$, with at least the specified 
ramification at the $P_i$.
\end{defn}

\section{Smoothing Families}\label{s-grd-families}

In this section we describe the families of curves whose limit linear
series we will study, called ``smoothing families'', and then give
some basic existence results. While the definition of a smoothing family is
rather technical, we expect that most applications will involve smoothing a
given reducible curve over a one-dimensional base, so we conclude with a
theorem giving the existence of such families satisfying all our technical
conditions, given the desired reducible fiber. However, we work over a fairly 
arbitrary base, because this allows the use of arguments in the universal
setting, in negative expected dimension, and in certain pathological cases 
where expected dimension is satisfied, but only over a positive-dimensional
``special fiber''.

Our central technical definition is:

\begin{defn}A morphism of schemes $\pi: X \rightarrow B$, together with 
sections $P_1, \dots,
P_n: B\rightarrow X$ constitutes a {\bf smoothing family} if:

\begin{Ilist}
\itm $B$ is regular and connected;
\itm $\pi$ is flat and proper;
\itm The fibers of $\pi$ are genus-$g$ curves of compact type;
\itm The images of the $P_i$ are disjoint and contained in the smooth
locus of $\pi$;
\itm Each connected component $\Delta'$ of the singular locus of $\pi$ maps
isomorphically onto its scheme-theoretic image $\Delta$ in $B$, and
furthermore $X|_{\pi ^{-1} \Delta}$ breaks into two (not necessarily
irreducible) components intersecting along $\Delta'$;
\itm Any point in the singular locus of $\pi$ which is smoothed 
in the generic fiber is regular in the total space of $X$;
\itm There exist sections $D_i$ contained in the smooth locus of 
$\pi$ such that every irreducible component of any geometric fiber of $\pi$ 
meets at least one of the $D_i$.
\end{Ilist}
\end{defn}

We begin with a lemma on two methods of obtaining new smoothing families
from a given one:

\begin{lem}Let $X/B$, $P_i$ be a smoothing family. Then
\begin{ilist}
\itm If $B' \rightarrow B$ is either a $k$-valued point of $B$ for any
field $k$, a localization of $B$, or a smooth morphism with $B'$ connected, 
then base change to $B'$ gives a new smoothing family.
\itm If $\Delta'$ is a node of
$X/B$ which is not smoothed in the generic fiber, let $Y, Z$ be the components
of $X$ with $Y \cup Z =X$, $Y \cap Z = \Delta'$. Then restriction to $Y$ or
$Z$ gives a new smoothing family.
\end{ilist}
\end{lem}

\begin{proof}For (i), the only properties of a smoothing family not preserved under arbitrary base change are (I) and (VI), which are easily checked in our specific cases. 

For (ii), the only condition that isn't immediately clear is that
flatness is preserved. However, this follows from the exact sequence of sheaves on $X$ 
$$0 \rightarrow \O_X \rightarrow \O_Y \oplus \O_Z \rightarrow \O_{Y \cap Z}
\rightarrow 0$$
together with the hypothesized flatness of $\O_X$ and $\O_{Y \cap Z}$ over $\O_B$.
\end{proof}

We now proceed to develop some results on construction of smoothing
families.

\begin{lem}\label{grd-make-fam} Let $\pi:X \rightarrow B$ be a family 
satisfying conditions
(I)-(III) and (VI) of a smoothing family, $\bar{X}_0$ a chosen geometric 
fiber of $\pi$ mapping to a point $P \in B$, and
$\bar{P}_i$ smooth closed points on $\bar{X}_0$ with images in $X_P$ having
residue fields separable 
extensions of $\kappa(P)$. Suppose further that each component $\Delta'$ of 
the singular locus of $X/B$ is flat over its image $\Delta$ in $B$. 
Then there is an 
\'etale base change of $\pi$ and sections $P_i$ specializing to the 
$\bar{P}_i$ which yield a smoothing family still containing $\bar{X}_0$ as a
geometric fiber, and with 
the same geometric generic fiber as $\pi$. 
\end{lem}

\begin{proof} We can Zariski localize $B$ to avoid any components of the
singular locus not occurring in $X_0$, and to insure that the sections we will contruct are disjoint and in the smooth locus. First, in addition to the $\bar{P}_i$, choose one 
smooth closed point $\bar{D}_i$ on each component of $X_0$, each having
field of definition a separable extension of $\kappa(P)$ (this is possible
by \cite[Prop. 2.2.16 and Cor. 2.2.13]{b-l-r}). Next, by 
\cite[Prop. 2.2.14]{b-l-r}, after possible \'etale base change we can 
find the desired sections $P_i$ and $D_i$ of $\pi$, each going through the 
corresponding $\bar{P}_i$ or $\bar{D}_i$. All that remains is to show that
we can obtain condition (V) as well. Since the singular locus of a family of
nodal curves is finite and unramified over the base, our flatness hypothesis
implies that each connected component $\Delta'$ is \'etale over its image
$\Delta$, and using \cite[Prop. 2.3.8 b)]{b-l-r} together with \cite[Cor.
V.1, p. 52]{ra1} in the case that $\Delta \neq B$, after an \'etale base change $\Delta'$ will map isomorphically to $\Delta$, giving the first half of (V).

Finally, we need to make sure that $X$ breaks
into components around each node. For each connected component $\Delta'$ of 
the singular locus of $\pi$, 
it suffices to produce an \'etale base change which causes the generic fiber
$X_1^{\Delta}$ of $X|_{\Delta}$
to break. By hypothesis, $X_1^{\Delta}$ breaks geometrically, and it will
break into components over an intermediate field $K$ if and only if 
the geometric components are $\Gal(\bar{K}/K)$-invariant. But this may be accompished after \'etale base change by \cite[Lem. 4.2]{dj-o1}, again using \cite[Cor. V.1, p. 52]{ra1} in the case $\Delta \neq B$.
\end{proof}

For typical applications of limit linear series, we expect that the
following theorem, which follows fairly easily from a theorem of Winters, 
will render irrelevant the technical hypotheses of our smoothing families:

\begin{thm}\label{grd-fam-exist}Let $X_0$ be any curve of compact type over an algebraically
closed field $k$, and $\bar{P}_1, \dots, \bar{P}_n$ distinct smooth closed 
points. 
Then $X_0$ may be placed into a smoothing family $X/B$ with sections $P_i$
specializing to the $\bar{P}_i$, where $B$ is a curve over $k$, and where 
the generic fiber of $X$ over $B$ is smooth.
\end{thm}

\begin{proof} Setting all $m_i=1$,
we can apply \cite[Prop. 4.2]{wi1} to obtain a proper map over 
$\Spec k$ from some regular surface $\tilde{X}$ to some regular curve 
$\tilde{B}$, having $X_0$ as a fiber. This must automatically be flat, and if we localize $\tilde{B}$ we can assume all fibers are at most nodal. We then claim that the
generic fiber $X_1$ must be smooth: indeed, all the local rings are regular by hypothesis, so by \cite[Cor. 16.21]{ei1} the residue fields of any non-smooth points would have to be inseparable over $K(B)$, which cannot happen in the case of nodes, since they are always unramified. Since the base field is algebraically closed, we need not worry about separability 
of residue field extensions on the closed fiber. Finally, our nodes are 
all isolated points, so the map to their image is a finite, unramified map of local schemes with algebraically closed residue field, and hence an isomorphism. Therefore, we can 
apply the preceding lemma to obtain our desired smoothing family.
\end{proof}

\begin{rem} There are a number of differences between our definition of a
smoothing family, and the one used in Eisenbud and Harris' original
construction in \cite{e-h1}. None of these are due to the different
construction. Extra conditions such as the reducedness of $B$ and the
regularity of $X$ at smoothed nodes are in fact necessary to ensure that
certain closed subschemes are actually Cartier divisors, and the condition 
that $X$ break into
distinct components above the nodes is likewise tacitly assumed, but not
automatic. The regularity of $B$ is necessary to make the sort of 
dimension-count arguments employed in the construction. Conversely, the 
hypotheses on the characteristic (or even existence) of a base field appears
to be unnecessary in their construction, as does the hypothesis that the
relatively ample divisor be disjoint from the ramification sections. The
only hypothesis we include here that may be truly gratuitous is that the
relatively ample divisor be composed of global sections, but it is
convenient and, as we have shown, not difficult to achieve. 
\end{rem}

\begin{rem}We do not claim that the moduli scheme could not be constructed
under weaker hypotheses, but merely that our hypotheses are those which
are necessary for our particular argument. It seems quite likely that one
could drop many of the hypotheses on both $X$ and $B$ if one carried out the construction in a universal setting and then pulled back the result to arbitrary families.
\end{rem}

\begin{rem}\label{grd-base-ex} It is not true that condition (VI) of a 
smoothing family is
preserved under base change by arbitrary closed immersions $B' \rightarrow
B$, even when $B'$ is regular and connected. For example, consider any
smoothing family with $B=\A^2_k$, and having a node $\Delta'$ with
$\Delta$ given by the $x$-axis. Then if $B'$ is the parabola $y=x^2$, base
change to $B'$ will create a singularity in $X$ above the origin.
\end{rem}

\begin{rem}In fact, the hypotheses for a smoothing family $\pi$ imply that 
every connected component of the singular locus of $\pi$ is regular, and in
particular irreducible and reduced. However, we will not need this, so we do not pursue it.
\end{rem}

\section{The Relative $\cG^r_d$ Functor}\label{s-grd-functor}

Given, in addition to a smoothing family, integers $r,d$, and ramification 
sequences $\alpha^i := \{\alpha_j^i\}_j$ for each of
our $P_i$, we will associate a $\cG^r_d$ functor to our smoothing family; 
this
functor will initially appear to include a lot of extraneous data, but we
will show that it actually gives the ``right'' functor, at least in the sense
that it associates a reasonable set to any geometric point of $B$. 

However, before defining the functor, we give some preliminary lemmas and definitions. In order to
ensure that our functor is globally well-defined, we will need the following easily-verified lemma.

\begin{lem}\label{grd-locokay}Let $\pi: X\rightarrow B$ be a proper morphism with 
geometrically reduced and
connected fibers, $\L$ and $\L'$ two isomorphic line bundles on $X$, and $V$ 
and $V'$ sub-modules of $\pi_* \L$ and $\pi_* \L'$ respectively. Then the 
property that ``$V$ maps into
$V'$'' is independent of the choice of isomorphism between $\L$ and $\L'$.
\end{lem}

We also describe a generalized notion of sub-bundle:

\begin{defn}\label{grd-sub-bundle} Let $\pi: X \rightarrow B$ be a morphism of schemes, and $\L$ a line bundle on $X$. A sub-sheaf $V$ is defined to 
be a {\bf sub-bundle} of $\pi_* \L$ if in addition to $V$ being a locally 
free sheaf, for any $S \rightarrow B$, the map 
$V_S \rightarrow \pi_{S*} \L_S$ remains injective. 
\end{defn}

Note that in this definition, we are pushing forward the pullback of $\L$, 
and not the other way around. The required sheaf map is gotten by composing 
the induced map $V_S
\rightarrow (\pi_* \L)_S$ with the natural map $(\pi_* \L)_S \rightarrow
\pi_{S*} \L_S$.

Finally, we define ramification conditions in this context.

\begin{defn} Let $X/B$ be a proper relative curve with line bundle $\L$ of degree $d$. Let $V$ be a sub-bundle of $\pi_* \L$ on $X/B$ of rank $r$, and $P$ a smooth section of $\pi$. Consider the sequence of maps
$$V \rightarrow \pi_* \L|_{(d+1)P} \rightarrow \pi_* \L|_{dP} \rightarrow
\dots, \pi_* \L|_{P} \rightarrow 0$$
We denote by $\beta_m$ the composition map $V \rightarrow \pi_* \L|_{mP}$.
Given a sequence of $r$ increasing integers $a_j$ between $0$ and $d$, 
we say that $V$ has {\bf vanishing sequence at least $\{a_j\}_j$}, or {\bf 
ramification sequence at least $\{a_j-j\}_j$}, if $\rk \beta_m \leq j$ for 
all $m \leq a^i_j$. 
\end{defn}

Finally, to simplify notation, and because it will be enough for inductive degenerations, we will restrict our families to reducible
curves with only two components. 

\begin{sit}\label{grd-sit}We assume that $X/B$ is a smoothing family with
at most one node (in the sense that 
the singular locus of $\pi$ is irreducible). If there is a node, we 
introduce some notation: denote by $\Delta'$ the singular locus of $\pi$, 
and $\Delta$ its image in $B$; by hypothesis, $\pi$ maps $\Delta'$ 
isomorphically to $\Delta$. We now distinguish three cases: case (1) is that
there is no node; case (2) is that $\Delta$ is all of $B$; and case (3) is
that $\Delta$ is a Cartier divisor on $B$. In cases (2) and (3), we denote by $Y$ and
$Z$ the components of $X|_{\pi^{-1}\Delta}$, necessarily smooth and intersecting along $\Delta'$.
\end{sit}

We observe that with the specified hypotheses, these three cases are all the
possibilities: indeed, since $B$ is regular, completing and examining the universal deformation of a nodal curve described \cite[p. 82]{d-m} easily shows that  if $\Delta$ is non-empty, it is 
locally generated principally. Note that with no hypotheses on the base, this would not be true; in fact, one can construct a families of nodal curves over a quadric cone having a node whose image is the union of three lines through the cone point.

In case (2), we will make use 
of the natural morphism (actually an isomorphism onto a connected component) 
$\Pic ^{d-i} (Y_T/T) \times \Pic^i(Z_T/T)\rightarrow \Pic^d(X_T/T)$ for any 
$i$ and any $T$ over $B$, in order to think of a pair of line bundles $\L_Y, \L_Z$ on $Y_T$
and $Z_T$ as a line bundle on $X_T$, which we will denote $(\L_Y, \L_Z)$. Note that this is only defined up agreement locally on the base, but this 
will not be a problem as we will consider the sheafified Picard functor.
In case (3), by the nonsingularity
hypothesis, $Y$ and $Z$ are Cartier divisors in $X$, so we have associated
line bundles on $X$, $\O_X(Y)$ and $\O_X(Z)$. Moreover, because $\Delta$ is a
Cartier divisor on $B$, and $\O_X(Y+Z) \cong \O_X(\pi ^* \Delta) \cong \pi ^* 
\O _B(\Delta)$, we have that locally on $B$, $\O_X(Y+Z) \cong \O_X$. 

Given a morphism $f: T \rightarrow B$, we denote by a subscript $T$ the various
pullbacks under $f$. We now describe our functor.

\begin{defn} The functor $\cG^r_d(X/B, \{(P_i, \alpha^i)\}_i)$ associates to $T/B$ the set of objects described as follows, modulo the equivalence induced by sheafification with respect to Zariski localization of $T$.
\begin{itemize}
\item[Case (1)] a line bundle $\L$ of degree $d$ on $X_T$, together with a 
rank $r$ sub-bundle $V$ of $\pi_{T*} \L$, having ramification sequence at
least $\alpha^i$ along the $P_{i,T}$.
\item[Case (2)] a line bundle $\L$ of degree $d$ on $X_T$, which has degree
$d$ when restricted to $Y_T$, and degree $0$ on $Z_T$, together with rank 
$r$ sub-bundles $V_0, \dots, V_d$ of $\pi_{T*} \L^i$, where $\L^i := 
(\L|_{Y_T} (-i\Delta'_T), \L|_{Z_T} (i \Delta'_T))$. Each $V_i$ must map 
to $V_{i+1}$ under the natural map given by inclusion on $Z_T$ and $0$ on 
$Y_T$, and each $V_i$ must map to $V_{i-1}$ under inclusion on $Y_T$ and 
$0$ on $Z_T$. Finally, we impose ramification along the $P_{i,T}$ as in case 
(1), with the caveat that we impose it only on $V_0$ if $P_i$ is on $Y$,
and only on $V_d$ if $P_i$ is on $Z$.
\item[Case (3)] a line bundle $\L$ of degree $d$ on $X_T$, which has degree $d$ when
restricted to $Y_T$, and degree $0$ on $Z_T$, together with rank $r$ sub-bundles $V_0, \dots, V_d$ of $\pi_{T*} (\L^i)$, where $\L^i := (\L \otimes \O_X(Y)_T^{\otimes i})$. Each $V_i$ must map to $V_{i+1}$ under the natural 
map $\pi_{T*} (\L^i) \rightarrow \pi_{T*} (\L^{i+1})$. Further, locally on $T$, we have $\O_X(Y+Z)_T \cong \O_{X,T}$, and we require that $V_i$ map to 
$V_{i-1}$ under the induced map. Finally, we impose
the desired ramification along the $P_{i,T}$ as in the first two cases, 
imposing it only on $V_0$ if $P_i$ specializes to $Y$, and only on $V_d$ 
if $P_i$ specializes to $Z$.
\end{itemize}
\end{defn}

\begin{rem} By Lemma \ref{grd-locokay}, our functor is well-defined despite the fact that we cannot distinguish line bundles on $X$ which are isomorphic locally on $B$. By the same token, the compatibility condition 
on $V_i$ in case (3) is independent of choice of local isomorphisms 
$\O_X(Y+Z) \cong \O_X$, and it isn't hard to see that the definition of 
$\cG^r_d$ in cases (2) and (3) is independent of the choice of $Y$ and $Z$.
\end{rem}

We also have:

\begin{defn} $\cG^{r, \sep}_d$ is the subfunctor of $\cG^r_d$ consisting of 
those linear series which are separable in every fiber. This is 
self-explanatory for smooth curves, while for reducible curves we require 
both $V_0|_Y$ and $V_d|_Z$ to be separable.
\end{defn}

One can verify quite directly that the $\cG^r_d$ we have defined is in
fact a functor. However, since we defined it differently in three separate 
cases, we also need to check:

\begin{lem}\label{grd-base-change}$\cG^r_d$ and $\cG^{r,\sep}_d$ are 
compatible 
with base change.
\end{lem}

\begin{proof} Suppose we pull back from $X/B$ to a new
smoothing family $X'/B'$. The only case with anything to check is if case
(3) pulls back to case (1) or (2). For the former, we note that $\O_X(Y)$
and $\O_X(Z)$ pull back to the trivial bundle, so all the maps between the
$\L^i$ are isomorphisms. For the latter, the point is that $\O_X(Y)$ clearly
pulls back to $\O_{Z'}(\Delta')$ on $Z'$, from which one checks that locally
on the base, $\O_X(Z)$ pulls back to $\O_{Z'}(-\Delta')$, and similarly with
$Y$ and $Z$ switched. 

Lastly, because $\cG^{r, \sep}_d$ was defined as a sub-functor of $\cG^r_d$ in
terms of behavior on fibers, it immediately follows that it too is
compatible with base change.
\end{proof}

\section{Representability}\label{s-grd-represent}

The main theorem is the representability of our $\cG^r_d$ functors, together with a lower bound on its dimension.
However, to ease the pain of the proof, we begin
with some technical lemmas before proceeding to the statement and proof of
the main theorem.

We begin with some compatibility checks on our notion of sub-bundle:

\begin{lem}\label{grd-sub}Our notion of sub-bundle has the following
desireable properties:
\begin{ilist}
\itm Suppose we have $\L$ such that $\pi_* \L$ is locally free, and the 
higher derived pushforward functors vanish. Then our definition of
sub-bundle of $\pi_* \L$ is equivalent to the usual one (that is, a locally
free sub-sheaf with locally free quotient).
\itm If $D$ is an effective Cartier divisor on $X$, flat over $B$, and $\L$
any line bundle on $X$, then a sub-sheaf $V$ of $\pi_* \L$ is a sub-bundle 
of $\pi_* \L$ if and only if it is a sub-bundle of $\pi_* \L(D)$ under the
natural inclusion.
\itm Let $V_1$, $V_2$ be sub-bundles of rank $r$ of $\pi_* \L$ in our
sense, and suppose $V_1 \subset V_2$. Then $V_1=V_2$.
\end{ilist}
\end{lem}

\begin{proof} For (i), we first note that by \cite[Cor. 6.9.9]{ega32} (see 
also \cite[6.2.1]{ega32}), the natural map 
$(\pi_* \L)_S \rightarrow \pi_* \L_S$ is an 
isomorphism. Now, if the quotient $Q := \pi_* \L/V$ is locally
free, we have $V_S \hookrightarrow (\pi_* \L)_S$ for any $S$ over $B$, 
and therefore $V_S \hookrightarrow \pi_* \L_S$, as desired.
Conversely, $V_S \hookrightarrow \pi_* \L_S$ for all $S$ means that the 
injectivity 
of $V \rightarrow \pi_* \L$ is preserved under base change; this in turn 
implies that $\Tor^1_S(\O_S, Q_S)=0$ for all $S$, since $\pi_* \L$ has 
vanishing $\Tor$. By \cite[Prop. 6.1]{ei1}, we conclude that $Q$ is flat, 
hence locally free, completing the proof of (i). 

Assertion (ii) will follow immediately if show that 
$\pi_* \L_S \rightarrow \pi_* \L(D)_S$ is injective for all $S$. However, by the flatness of $D$ over $B$, the cokernel of 
$\L \hookrightarrow \L(D)$ is flat over $B$, so injectivity of this map is preserved under base change, and applying $\pi_*$ gives the desired result.

Finally, (iii) is straightforward: let $Q = V_2/V_1$, and let $b \in B$ be 
any point of $B$. If we base change to $\Spec \kappa(b)$, we get find from the definition of sub-bundle that $V_{1b} \hookrightarrow V_{2b}$, so since their dimension is the same, we get $Q_b=0$, and by Nakayama's lemma we conclude $Q=0$ and $V_1=V_2$, as asserted.
\end{proof}

We also have a lemma illustrating how we will use our sections $D_i$:

\begin{lem}\label{grd-ample}Let $X/B$ be a smoothing family, and $\L^i$ any
finite collection of line bundles on $X$,
of degree $d$. Then there exists an effective divisor $D$ on 
$X$ satisfying:
\begin{ilist}
\itm $D$ is flat over $B$, and supported in the smooth locus of $\pi$.
\itm $\pi_* \L^i(D)$ is locally free and $R^1 \pi_* \L^i(D)=0$ for all $i$.
\itm $\pi_* \L^i(D) \rightarrow \pi_* \L^i(D)|_{mP_j}$ is surjective for
all $i,j$ and all $m \leq d+1$.
\itm In case (3) of Situation \ref{grd-sit}, $D$ may be written as $D^Y+D^Z$, with $D^Y|_{\pi^{-1}\Delta}$ contained in $Y$, and similarly for $D^Z$.
\end{ilist}

Furthermore, \'etale locally on $B$ we may require that $D$ is disjoint from the $P_i$ as well.
\end{lem}

\begin{proof}
With $D_i$ any collection of sections as in the definition of a smoothing 
family, let $D' = \sum _i D_i$; then $D'$ is $\pi$-ample, so locally on 
$B$, for $\ell$ sufficiently large, $\ell D'$ will have the desired properties. Since $B$ is Noetherian, we can choose such an $\ell$ globally.

The \'etale-local disjointness assertion is obtained by constructing new sections \'etale locally as in the proof of Lemma \ref{grd-make-fam}, choosing the $\bar{D}_i$ to be distinct from the $\bar{P}_i$.
\end{proof}

We can now prove our central result:

\begin{thm}\label{grd-main}If $\pi:X \rightarrow B$, $P_1, \dots, P_n: 
B \rightarrow X$ is a
smoothing family satisfying the two-component hypothesis of Situation
\ref{grd-sit}, and $\alpha^i:=\{\alpha^i_j\}_j$ ramification sequences, 
then $\cG^r_d=\cG^r_d(X/B;
\{(P_i, \alpha^i)\}_i)$ is represented by a scheme $G^r_d$, compatible with 
base change to any other smoothing family.
This scheme is projective, and if it is non-empty, the local ring at
any point $x \in G^r_d$ closed in its fiber over $b \in B$ has dimension at 
least $\dim \O_{B,b} + \rho$, where $\rho=\rho(g,r,d; \alpha^i)$ as in
Definition \ref{grd-expected}. Furthermore, $\cG^{r, \sep}_d$ is also 
representable, and is naturally an open subscheme of $G^r_d$.
\end{thm}

\begin{proof}Once the $\cG^r_d$ functor has been defined, the proof of its
representability is long but for the most part extremely straightforward, 
using nothing more than the well-known representability of the various
functors in terms of which we have described $\cG^r_d$. The one trick, borrowed from Eisenbud and Harris, is to twist a universal line bundle 
$\L$ by a high power of an ample divisor so sub-bundles of its pushforward are parametrized by a standard Grassmannian
scheme. The dimension count is an altogether different story; it is 
harder than in Eisenbud and Harris' construction, and is essentially the
subject of the final section this chapter. 

We note that our functor is visibly a Zariski sheaf, so we can check representability Zariski locally on $B$. Furthermore, it is clear that imposing ramification conditions, consisting of imposing rank conditions on sequences of maps of locally free sheaves, give closed subfunctors, so it suffices to check representability without imposing ramification. 

As in defining the functor, we have three cases to consider. The first is
the simplest. We start in this case with the relative Picard scheme
$P=\Pic^d(X/B)$, obtained for instance from \cite[Thm. 9.4.1]{b-l-r} and
twisting by sections to obtain the desired degree. We denote by $\tilde{\L}$
the universal line bundle on $X \times_B P$. Let $D$ be the divisor provided by Lemma \ref{grd-ample} for $\tilde{\L}$, 
viewing $X \times_B P$ as a smoothing family over $P$. Now, we let $G$ be 
the relative Grassmannian scheme of $\pi_{P*} (\tilde{\L}(D))$. 
We define our $G^r_d$ scheme to be the closed subscheme
of $G$ cut out by the condition that any sub-bundle $V$ of
$\pi_{P*} (\tilde{\L}(D))$ vanishes on $D$, which naturally gives a closed subscheme cut out locally by
minors. This completes the construction in the first case.

Now, in the second case, we use the Picard schemes $P^i := \Pic
^{d-i,i}(X/B)$, the schemes parametrizing line 
bundles on $X$ with degrees $d-i$ and $i$ when restricted to $Y$ and $Z$ 
respectively. These are all naturally isomorphic to one another by twisting the line bundles on each component by $\Delta'$ and $-\Delta'$; in
particular, we can identify all of them with a fixed $P$ over $B$. 
On each $P^i$, we have a universal line
bundle $\tilde{\L}^i$, and just as in the first case, we take a very ample
divisor $D$ obtained from Lemma \ref{grd-ample} for the $\tilde{\L}^i$,
twist $\tilde{\L}^i$ by $D$, and then construct Grassmannian bundles
$G^i$, this time one for each $\tilde{\L}^i$. Denoting by $G$ the product of 
all these Grassmannians over $P$, we take the closed subscheme inside $G$ cut 
out by, as in the first case, vanishing on $D$ and the required
ramification conditions along the $P_i$. Here, we actually write
$D=D^Y+D^Z$, where $D^Y$ and $D^Z$ are supported on $Y$ and $Z$
respectively, and impose vanishing along $D^Y$ only in $G^0$, and along 
$D^Z$ only in $G^d$. Finally, we make use of the construction Lemma 
\ref{grd-pre-lg} to
add the requirement that the $V_i$ each map into $V_{i+1}$ on $Z$ and 
$V_{i-1}$ on $Y$ under the natural maps, also as in the definition of the
functor. This completes the construction in the second case.

In the final case, the first step is to work sufficiently locally on $B$
that $\Delta$ is principal, so that $\O_B(\Delta) \cong \O_B$ and $\O_X(Y+Z)
\cong \O_X$, and fix a choice of this isomorphism. The rest proceeds very similarly to the second
case: our Picard schemes $P^i$ are described identically, but now to
describe isomorphisms between the $P^i$, we tensor as necessary by 
$\O_X(Y)$. Replacing the maps between the
$\tilde{\L}^i$ with the appropriate maps for this case, the rest of the 
construction then proceeds identically to the previous case.

Because in each case the construction used only Picard schemes,
Grassmannians, fiber products, and closed subschemes obtained by bounding
the rank of maps between vector bundles, it nearly follows from the standard
representability theorems for these functors that the
$G^r_d$ scheme we have constructed represents the $\cG^r_d$ functor. We do
need to note that in the second and third cases, our conditions for vanishing 
along $D$ actually imply that all $V_i$ vanish along $D$: in the second
case, this follows simply because $D^Y$ and $D^Z$ are disjoint from $\Delta'$; in the third case, we have similarly that $D^Y$ is disjoint from $Z$ and $D^Z$ disjoint from $Y$. By Lemma \ref{grd-sub}, our definition of sub-bundle is
compatible with the usual definition for the Grassmannian functor. We have thus proven representability. It easily follows from the projectivity of Grassmannians and Picard schemes that our $G^r_d$ scheme is projective over $B$. Lastly, compatibility with base 
change has already been proven in Lemma \ref{grd-base-change}.

We now verify that the moduli scheme we have constructed has the desired
lower bound on its dimension. Our ambient scheme $G$ is a product of 
Grassmannians over a Picard scheme, so since $B$ was assumed to be regular, we conclude that $G$ is regular. Hence, in order to bound the codimension 
of $G^r_d$ in $G$ it suffices to consider the codimensions of 
each condition cutting it out. We denote by $d''$ the rank of $\pi_{P*} \tilde{\L}^i(D)$; we have $d''=d+\deg D + 1-g$.

Now, in the first case, vanishing along $D$ imposes $(\rk V)(\rk
\pi_{P*} \tilde{\L}|_D)=(r+1)(\deg D)$ conditions. Next, we consider the codimension of the ramification conditions. For this calculation, it suffices to work \'etale locally, so by Lemma \ref{grd-ample} we may assume that $D$ is disjoint from the $P_i$, in which case it follows that ramification conditions imposed
on sub-bundles of $\L(D)$ are equivalent to the desired ramification for 
sub-bundles of $\L$. 
Since the evaluation maps
$\pi _{P*} \tilde{\L}(D) \rightarrow \pi_{P*} (\tilde{\L}(D)|_{jP_i})$ are
surjective, each condition defines a Schubert cycle. In particular,
by \cite[Thm. 6.3, Cor. 5.12 (b)]{b-v} the imposition of ramification at $P_i$ gives an
integral subscheme of codimension $\sum_j (a^i_j-j)=\sum_j (\alpha^i_j)$
inside $G$. Thus the total codimension of any component of $G^r_d$ inside $G$ 
is at most $(r+1)(\deg D)+ \sum_j(\alpha^i_j)$. 

In the second and third cases, the only real difference is that we replace 
the Grassmannian with the linked Grassmannian of Appendix \ref{s-grd-lg}; it is
easily verified that because the maps on $\pi_{P*} \L^i$ are induced from maps
on the $\L^i$, they satisfy the conditions of a linked Grassmannian
(Definition \ref{grd-lg-def}). Note that in the case of a reducible fiber, 
everything in the kernel of $f_i$ really is in the image of $g_i$ and vice versa, because the $\L^i$ have all been constructed to be sufficiently ample.
Then it follows from Theorem \ref{grd-main-lg} that every component of the
linked Grassmannian has codimension $d(r+1)(d''-r-1)$, and
the rest of the calculation proceeds the same way, with the minor exception
that we have to compute vanishing on $D^Y$ and $D^Z$
separately and use $\deg D^Y + \deg D^Z = \deg D$.

Finally, we conclude the desired dimension statement, noting that $G$ is 
catenary, and smooth over $B$ of relative dimension 
$(r+1)(d''-r-1)+g$ in the first case and $(d+1)(r+1)(d''-r-1)+g$ in the 
second and third cases.

Lastly, we need to show that the sub-functor of separable limit series is
representable by an open subscheme. Denote by $\tilde{\F}^i$ the universal sub-bundles of 
$\pi_{P*} \tilde{\L}$ on our $G^r_d$ scheme. As in the proof of Proposition 
\ref{grd-plucker}, we can construct a map $\tilde{\F}^i \otimes \O_{X \times _B
G^r_d} \rightarrow \sP^r(\tilde{\L})$ where $\sP^r$ denotes the bundle of
principal parts of order $r$; taking $(r+1)$st exterior powers gives a map
$s^{\univ}: \det (\tilde{\F}^i) \rightarrow \tilde{\L}^{\otimes r+1}\otimes
(\Omega^1_{X/B})^{\otimes \binom{r+1}{2}}$; we already noted that in the smooth 
case, our separable subscheme is the image under $\pi_P$ in $G^r_d$ of the complement of the closed subscheme cut out as the kernel of $s^{\univ}$. For the second and third cases, we restrict to the smooth locus to avoid the problem that $\Omega^1_{X/B}$ is no longer locally free. It is then easy to check that the same construction applied to $\F^0$ and $\F^d$ will give the subscheme of linear series which are separable on $Y$ and $Z$ respectively, and their intersection gives the desired $G^{r, \sep}_d$ subscheme.
\end{proof}

Our first application is the same regeneration/smoothing theorem due to
Eisenbud-Harris, except that now it {\it a priori} gives results on 
smoothings of
crude limit series as well, and we are also able to include upper bounds of
dimensions of general fibers in certain cases. 

We have:

\begin{cor}\label{grd-main-cor}In the situation of Theorem \ref{grd-main}, 
suppose that $\rho
\geq 0$, that $U$ is
any open subscheme of our $G^r_d$ scheme, and that for some point $b \in B$, 
the fiber of $U$ over $b$ has the expected dimension $\rho$. Then every
point of the fiber may be smoothed to nearby points. Specifically:
\begin{ilist}
\itm The map from $U$ to $B$ is open at any point in the fiber over
$b$, and for any component $Z$ of $U$ whose image contains $b$, the generic 
fiber of $Z$ over $B$ has dimension $\rho$. 
\itm If further $U$ is closed in $G^r_d$, then there is a neighborhood
$V$ of $b$ such that the preimage of $V$ in $U$ is open over $V$, and for 
each component $Z$ of $U$, every component of 
every fiber of $Z$ over $V$ has dimension precisely $\rho$.
\end{ilist}

In particular, if $X_0$ is a curve of compact type (with two components)
over an algebraically closed field, with $\bar{P}_1, \dots, \bar{P}_n$ 
distinct smooth closed points of $X_0$, $\alpha^i$ any collection of
ramification sequences, and $U_0$ any open subset of
$G^r_d(X_0/k; \{(P_i, \alpha^i)\}_i)$ having expected dimension $\rho$, then
there exists a smooth curve $X_1$ over a one-dimensional function field
$k'$ over $k$, specializing to $X_0$, with points $P_i$ specializing to the
$\bar{P}_i$, and such that every point of $U_0$ smooths 
to $X_1$; if further $U_0=G^r_d(X_0/k; \{(\bar{P}_i, \alpha^i)\}_i)$, then 
$G^r_d(X_1/k'; \{(P_i, \alpha^i)\}_i)$ also has dimension $\rho$.
\end{cor}

\begin{proof} For (i), let $x \in Z$ be any
closed point in the fiber of $Z$ over $b$, and $\eta$ the generic point of
$Z$. Say $\eta$ maps to $\xi$; then the dimension of the fiber of $Z$ over
$\xi$ is at most $\rho$, by \cite[Thm. 13.1.3]{ega43}, and at least $\rho$ by Theorem
\ref{grd-main} after base change to $\xi$. But because the fiber over $b$ has dimension $\rho$, and by Theorem \ref{grd-main} we have $\dim
\O_{Z,x} \geq \dim \O_{B,b} + \rho$, the image of $Z$ cannot have dimension less than $\dim \O_{B,b}$, so $\xi$ is the generic point of $B$.

For the openness assertion, it suffices to prove that the image of $U$ 
contains a neighborhood of $b$, since if we replace $U$ by any neighborhood
of a point of the fiber of $U$ over $b$, the hypotheses of our corollary are
still satisfied. Let $b_1$ be a point of $B$, specializing to $b$; let $B_1$ 
be the closure of $b_1$ in $B$, and consider the base change 
$U_1 \rightarrow B_1$. If $B_1$ has codimension $c$ in $B$, then every 
component of $U_1$ would have codimension at most $c$ in $U$, so if we restrict to a
component $Z_1$ of $U_1$ passing through $x$, we have $\dim \O_{Z_1,x} \geq 
\dim \O_{B_1,b} + \rho$, so arguing as before we see that $b_1$ must be in
the image of $U$. Now, by constructibility of the image, $f(U)$ must
contain some neighborhood of $b$, as desired. 

For (ii), if $U$ is closed in $G^r_d$ we have that it is proper over $B$, and
every component $Z$ of $U$ either contains $b$ in its image, or is supported
on a closed subset of $B$ away from $b$. If $Z$ maps to $b$, we can apply (i)
to conclude that $Z$ maps surjectively to $B$, and by \cite[Cor.
13.1.5]{ega43} the locus on $B$ of fibers of $Z$ having a
component of dimension greater than $\rho$ is closed, so taking its 
complement and intersecting over the finitely many components of $U$ gives 
a $V$ of the desired form. Openness then follows from (i) and the fact that 
all fibers over $V$ have dimension $\rho$.

Finally, given an $X_0$ as described, we can apply Theorem
\ref{grd-fam-exist} to place $X_0$ into a smoothing family $X/B$ with
generic fiber $X_1$; the desired assertions then follow immediately from the
main assertions of the corollary. 
\end{proof}

The finite case is particularly nice, but we put off any discussion of it
until after we have introduced the language of Eisenbud-Harris limit series 
in the next section.

Even without knowing anything about the separable locus being closed, which
in general seems to be a subtle issue, we can still obtain results on 
lifting from characteristic $p$ to characteristic $0$. However, note that the 
expected dimension hypothesis in the following corollary is not only key to 
the argument, but at least in some cases both non-vacuous and necessary for
the validity of the conclusion. See in particular 
\cite[Prop. 5.4, Rem. 8.3]{os7}.
In any case, our machinery now easily yields:

\begin{cor}\label{grd-lift} In the situation of Theorem \ref{grd-main}, suppose 
that $\rho \geq 0$, that $B$ is a mixed-characteristic DVR, and that the 
special fiber of some $U$ open inside $G^r_d$ has the expected 
dimension $\rho$. Then every point 
$x_0$ of $U$ in the special fiber may be lifted to characteristic 
$0$, in the sense that there will be a point $x_1$ of the generic fiber
of $U$ (and in particular of $G^r_d$) specializing to $x_0$. 

In particular, suppose that $X_0$ is a smooth, proper
curve over a perfect field $k$ of characteristic $p$, with
$\bar{P}_1, \dots, \bar{P}_n$ distinct closed 
points of $X_0$, $\alpha^i$ any collection of
ramification sequences, and $U_0$ any open subset of
$G^r_d(X_0/k; \{(\bar{P}_i, \alpha^i)\}_i)$ having expected dimension $\rho$;
then there exists a smooth curve $X_1$ over the fraction field $K$ of the Witt
vectors of $k$, specializing to $X_0$, with points $P_i$ specializing to the
$\bar{P}_i$, and such that every point of $U_0$ may be lifted to a point of 
$G^r_d(X_1/K; \{(P_i, \alpha^i)\}_i)$.
\end{cor}

\begin{proof}
The first assertion follows immediately from the openness proven in
Corollary \ref{grd-main-cor}.

For the second assertion, let $A$ be the Witt vectors of $k$; then by 
\cite[11, Thm 1.1]{b-l-r2} we can find an $X$ over $\Spec A$ whose special 
fiber is $X_0$, and since $A$ is complete and the $\bar{P}_i$ are smooth 
points we can lift them to sections $P_i$ of $X$. Applying the first 
assertion then gives the desired result.
\end{proof}

\begin{rem} This last corollary, dealing only with smooth curves, has nothing to do with limit linear series, and only uses the elementary lower bound on the dimension of a standard $G^r_d$ space with imposed ramification.
\end{rem}

\section{Comparison to Eisenbud-Harris Theory}\label{s-grd-compare}

This is all well and good, but our description of the limit series
associated to a reducible curve is rather cumbersome, so we now establish
the relationship to Eisenbud and Harris' limit series in this situation.

\begin{sit}\label{grd-eh-sit} $X/B$ is a smoothing family with $X$ 
reducible; specifically, falling into case (2) of Situation \ref{grd-sit}.
\end{sit}

The first step is to consider the ``forgetful'' map from our 
$\cG^r_d(X/B)$ functor into the product of $\cG^r_d(Y/B)$ and 
$\cG^r_d(Z/B)$.

We have:

\begin{lem}\label{grd-fr} In Situation \ref{grd-eh-sit}, given any $T$-valued point $\{(\L^i, V_i)\}_i$ of
$G^r_d(X/B)$, the pair $((\L^0|_{Y_T}, V_0|_{Y_T}), 
(\L^d|_{Z_T}, V_d|_{Z_T}))$ gives a $T$-valued point of 
$G^r_d(Y/B) \times G^r_d(Z/B)$. In particular,
this defines a morphism 
$FR: G^r_d(X/B) \rightarrow G^r_d(Y/B) \times _B G^r_d(Z/B)$. A limit series 
in $G^r_d(X/B)$ is separable if and only if its image under $FR$ is 
separable in both $G^r_d(Y/B)$ and $G^r_d(Z/B)$.
\end{lem}

\begin{proof}It is clearly enough to show that given $T/B$, and 
$\{(\L^i, V_i)\}_i$ a $T$-valued point of $G^r_d(X/B)$, then
$V_0|_{Y_T}$ is a sub-bundle of $\L^0|_{Y_T}$ and correspondingly for
$Z_T$ and $V_d$. This follows immediately from our definition of sub-bundle, since if a section of $V_{0S}$ vanishes on $Y_S$ for any $S$ over $T$, it defines a section of the negative line bundle $\L^0|_Z(-\Delta')$ and hence vanishes on $Z_S$ as well, and similarly for $V_{dS}$ with $Y$ and $Z$ switched.

The statement on separability is immediate from the definition of
separability of a limit series on a reducible curve.
\end{proof}

\begin{notn} In the same situation as the previous lemma, we denote the image under $FR$ of $\{(\L^i, V_i)\}_i$ by $((\L^Y, V^Y),(\L^Z, V^Z))$.
Further, we denote by
$V^Y_i$ the image of $V_i$ inside $\pi_*(\L^i|_Y) \cong 
\pi_*(\L^Y(-i\Delta'))$, and similarly for $Z$.
\end{notn}

\begin{lem}\label{grd-random-easy}In the same situation as the previous lemma, 
we have the following additional observations (and consequent notation):
\begin{ilist}
\itm $V^Y_i$ injects naturally into $V^Y$, and similarly for
$Z$;
\itm $V^Y_i$ will be contained in $\ker \beta^Y_i \subset V^Y$,
where $\beta^Y_i: V^Y \rightarrow \pi_{T*} \L^Y|_{i\Delta'}$ is the natural
$i$th order evaluation map at $\Delta'$, and $V^Z_i$ will 
similarly be contained in $\ker \beta^Z_{d-i} \subset V^Z$;
\itm The induced map $V_i \rightarrow V^Y \oplus V^Z$ in fact
exhibits $V_i$ as a sub-bundle (in the usual sense) of $V^Y \oplus V^Z$.
\end{ilist}
\end{lem}

\begin{proof}Assertions (i) and (ii) are clear. For (iii), it suffices to show that
$V_i \rightarrow V^Y \oplus V^Z$ is injective after any base change $S
\rightarrow T$, and this is easily verified from the definitions.
\end{proof}

\begin{defn}In case (2) of Situation \ref{grd-sit}, we define an 
{\bf Eisenbud-Harris (crude) limit series} on $X$ to be a 
pair $((\L^Y, V^Y),(\L^Z, V^Z))$ in $G^r_d(Y) \times _B G^r_d(Z)$ satisfying
$a^Y_i(\Delta')+a^Z_{r-i}(\Delta')\geq d$ for all $i$ (see below). The 
closed subscheme
of $G^r_d(Y) \times _B G^r_d(Z)$ obtained by these ramifications conditions
will be denoted 
$G^{r}_{d, \EH}(X/B)$. We also define $G^{r, \sep}_{d,\EH}(X/B) \subset
G^{r}_{d, \EH}(X/B)$ to be the 
open subscheme of limit series which are separable on each component, and 
$G^{r, \refn}_{d, \EH}(X)$ to be the open subscheme of {\bf refined} 
Eisenbud-Harris 
limit series satisfying $a^Y_i(\Delta') +a^Z_{r-i}(\Delta') = d$ for all 
$i$, or more
precisely, the complement of the closed subscheme satisfying 
$a^Y_i(\Delta') +a^Z_{r-i}(\Delta') > d$ for some $i$.
\end{defn}

We remark that these ramification conditions do in fact
give a canonical closed subscheme structure:
for each sequence of $r+1$ non-decreasing integers $0 \leq a_i \leq d$, we
get a closed subscheme defined by the conditions
$a^Y_i(\Delta') \geq a_i$, $a^Z_{r-i}(\Delta') \geq d-a_i$; there
are only finitely many such sequences, so the union of the closed subschemes
obtained over each of them is again a closed subscheme. However, this
definition gives us trouble when we attempt to show that our
$G^r_d(X/B)$ maps into $G^r_{d, \EH}(X/B)$, as it is difficult to describe
the $T$-valued points of a union of schemes in terms of the $T$-valued
points of the individual schemes. As a result, we settle for the following slightly weaker statement.

\begin{prop}\label{grd-eh1} We have the following facts about the image of
$FR: G^r_d(X/B) \rightarrow G^r_d(Y/B) \times _B G^r_d(Z/B)$:
\begin{ilist} 
\itm $FR$ has set-theoretic image precisely $G^{r}_{d, \EH}(X/B)$; 
\itm Scheme-theoretically, $G^r_d(X/B)$ maps into the closed subscheme 
satisfying for all $j$
$$a^Y_j(\Delta') + a^Z_{r-j}(\Delta')\geq d-1;$$ 
\itm The open subscheme of $G^r_d(X/B)$ mapping set-theoretically into 
$G^{r, \refn}_{d, \EH}(X)$ actually maps scheme-theoretically into 
$G^{r, \refn}_{d, \EH}(X) \subset G^{r}_{d, \EH}(X/B)$. 
\end{ilist}
\end{prop}

\begin{proof}In general, for a $T$-valued pair $((\L^Y, V^Y),(\L^Z, V^Z))$, 
define $a_j^Y$ to be the largest 
integer $i$ with $\rk \beta^Y_i \leq j$ everywhere on $T$, and similarly
for $Z$. The set-theoretic statement may be checked point by point, and is
equivalent to saying that when $T = \Spec k$ for some $k$,
$((\L^Y,V^Y),(\L^Z, V^Z))$ is
in the image of $FR$ if and only if $a_j^Y + a_{r-j}^Z \geq d$ for all $j$.
For (ii), it is enough to check that for arbitrary local
$T$, $a_j^Y+a_{r-j}^Z \geq
d-1$ for all $j$, and for (iii), we want to show in this case that if the
point obtained by restriction to the closed point of $T$ satisfies 
$a_j^Y+a_{r-j}^Z=d$ for all $j$, then the entire $T$-valued point does. 
In all cases, we make use of the fact from Lemma 
\ref{grd-random-easy} that $V^Y_i$ may be considered as lying inside 
$\ker \beta^Y_i$, and $V^Z_i$ in $\ker \beta^Z_{d-i}$. Conceptually, 
the basic idea is
that for $V_i$ to maintain rank $r+1$ at each $i$, the ranks of $\ker
\beta^Y_i$ and $\ker \beta^Z_{d-i}$ must add up to at least $r+1$, and
looking at $i = a^Y_j$ for different $j$ should yield the desired
inequalities. As we will see, this works over a field, but is not quite so
nice for a more general $T$.

Now, for the set-theoretic statement (i), suppose we have a
$k$-valued point $\{(\L^i,V_i)\}_i$ of $G^r_d(X/B)$; we
first show that it maps into $G^r_{d,\EH}(X/B)$. Since $V_i$ is glued 
from subspaces of $\ker \beta^Y_i$ and $\ker \beta^Z_{d-i}$
and has dimension $r+1$, we conclude that $\dim \ker \beta^Y_i +
\dim \ker \beta^Z_{d-i} \geq r+1$, 
so $\rk \beta^Y_i + \rk \beta^Z_{d-i} \leq r+1$, and 
it follows that $a_j^Y + a_{r+1-j}^Z
\geq d$ for all $j$. On the other hand, for a given $j$, set $i = a_j^Y$; we 
know that $\rk \beta^Y_{i+1} > j$, so one of the sections in $V^Y_i$ is
non-vanishing at $\Delta'$ when considered as a section of $\L^i(-i\Delta')$, 
and to use it in $V_i$, it must be glued to a section of $V^Z_i$
similarly non-vanishing at $\Delta'$. Thus, $\dim V_i < \dim \ker \beta^Y_i + \dim \ker \beta^Z_{d-i}$, so our earlier argument
gives $\rk \beta^Y_i + \rk \beta^Z_{d-i} < r+1$, hence 
$a_j^Y + a_{r-j}^Z \geq d$. For later use, note that when 
$a^Y_j+a^Z_{r-j}=d$
for all $j$, this argument shows that we have $\ker \beta^Y_i = V^Y_i$ and
$\ker \beta^Z_{d-i} = V^Z_i$ for all $i$, and in particular when $i=a^Y_j =
d-a^Z_{r-j}$, we get $\dim V^Y_i + \dim V^Z_i = r+2$.

Conversely, given a $((\L^Y,
V^Y),(\L^Z, V^Z))$ satisfying the Eisenbud-Harris inequalities, we 
construct the $\L^i$ by
gluing $\L^Y (-i\Delta')$ and $\L^Z((i-d)\Delta')$, and set $V_0 = V^Y$, 
$V_d = V^Z$. Note that sections in $V^Y$ which vanish at $\Delta'$ are
extended by $0$ along $Z$. If there is a non-vanishing section, then 
$a^Y_0=0$, so $a^Z_r\geq d$, and $\L^Z$ is necessarily 
$\O_{Z_T}(d\Delta')$, so $\L^i|_Z \cong \O_{Z_T}$, and we can (uniquely) 
extend sections not vanishing at $\Delta'$, also. We also observe that this 
implies that we have $V_0$ mapping into $V^Z$ under iterations of $f_i$. By 
symmetry, we can make the same arguments for $V^Z$ to get our $V_d$.
Now we inductively construct each $V_i$ for $i=1, 2, \dots, d-1$ in 
terms of $V_{i-1}$ and $V_d$. Our induction hypothesis will be that $V_{i-1}$ is linked to the previous $V_j$ and $V_d$ under iterates of $f_j$ and $g_j$, and furthermore that each $V_j$ has a basis of sections each 
of which is either non-vanishing at $\Delta'$, or vanishes uniformly on 
either $Y$ or $Z$, with at most one basis element in the first category. We denote the number of each of 
these by $r_j^1$, $r_j^2$, and $r_j^3$ respectively, where we have $r_j^1$ 
always $0$ or $1$, and $r_j^1+r_j^2+r_j^3 = r+1$ for all $j$. Finally, we also impose in our induction hypothesis that $r_j^3$ is always the maximal 
possible value, which is $\dim \ker \beta^Y_{j+1}$. Note that since this is
non-decreasing, if we construct a $V_i$ with $r_i^3=r_{i-1}^3$, maximality
is automatically satisfied. 

Now, for general $i$, suppose we have constructed the $V_j$ up to $V_{i-1}$
satisfying our induction hypothesis. To construct $V_i$, the basis elements 
vanishing on $Y$ must contain $f_{i-1} (V_{i-1})$, which is an
$(r_{i-1}^1 +r_{i-1}^2)$-dimensional space, and of course they must map into
$V^Z$. Since $f_{i-1} (V_{i-1})$ maps into $V^Z$, we can choose 
$r_{i-1}^1+r_{i-1}^2$ such sections, by taking any basis of $f_{i-1}
(V_{i-1})$. Next, the basis elements vanishing on $Z$ must be contained in 
$g_{i-1}^{-1} (V_{i-1})$, and we choose them to be a basis of the 
subspace of $g_{i-1}^{-1} (V_{i-1})$ vanishing on $Z$. 
This is at most an $r_{i-1}^3$-dimensional space, with equality if all of 
the $r_{i-1}^3$ basis elements of $V_{i-1}$ vanish to order greater than one 
at $\Delta'$. If there was a section vanishing 
to order exactly one at $\Delta'$, $g_{i-1}^{-1}$ will instead be 
$(r_{i-1}^3-1)$-dimensional. Now, by our induction hypothesis, the 
iterated image of $V_d$ under the $g_j$ is contained in $V_{i-1}$,
necessarily in the span of the basis elements vanishing on $Z$. Moreover,
since $i<d$, this image lies in the subspace of $V_{i-1}$ vanishing 
to order at least $2$ at $\Delta'$, so it 
is automatically contained in the span of the basis elements we have 
chosen for $V_i$ which vanish on $Z$. Now, if we had $r_{i-1}^3$ such 
basis elements, we are done. If not, we had a section of $V_{i-1}$ vanishing 
on $Z$ and vanishing to first order at $\Delta'$ on $Y$, so it follows that 
$\dim \ker \beta^Y_i = \dim \ker \beta^Y_{i+1} + 1$, and therefore that
$a^Y_{r+1-r^3_{i-1}}=i$; in particular, the required maximality of $r^3_i =
r^3_{i-1}-1$
is satisfied. It also follows that $a^Z_{r^3_{i-1}-1} \geq d-i$; if it is 
equal, we can find a section of $V^Z$ vanishing to order precisely $d-i$ at
$\Delta'$, which we could glue to our final section of $g_{i-1}^{-1} 
(V_{i-1})$ to obtain our $(r+1)$st generator for $V_i$, which will be 
non-vanishing at $\Delta'$. Otherwise, we have $a^Z_{r^3_{i-1}-1} > d-i$, 
so following through the definitions, $\dim \ker \beta^Z_{d-i+1} \geq 
r+2-r^3_{i-1} = 1+r^1_{i-1}+r^2_{i-1}$, and we can choose an 
$(r^1_{i-1}+r^2_{i-1}+1)$st generator vanishing on $Y$ to be our $(r+1)$st 
generator for $V_i$. This completes the proof of the set-theoretic 
surjectivity of $FR$ onto $G^{r}_{d, \EH}(X/B)$.

For the scheme-theoretic statements, let $T = \Spec A$ where $A$ is any local
ring with maximal ideal $\m$. By Lemma \ref{grd-random-easy}, for each 
$i$ we have 
$0 \rightarrow V_i \rightarrow V^Y \oplus V^Z \rightarrow Q \rightarrow 0$ 
for some free $Q$. Working modulo $\m$ and 
then using Nakayama's lemma, we find that for some $j$ depending on $i$, we 
can construct 
sub-bundles of $V_i$ whose images in $V^Y_i$ and $V^Z_i$ are sub-bundles of 
$V^Y$ and $V^Z$ of rank $r+1-j$ and $j$. These are contained in $\ker \beta^Y_i$ and $\ker \beta^Z_{d-i}$, and
we conclude that $a^Y_j + a^Z_{r+1-j} \geq d$. As in the fields case, if
we set $i = a^Y_j$, then by hypothesis $\rk \beta^Y_{i+1}$ is not
less than or equal to $j$ on all of $T$, so our constructed sub-bundle of
$V^Y$ for $i+1$ could have rank at most $r-j$, and the sub-bundle of $V^Z$ 
would have
to have rank at least $j+1$, giving $\rk \beta^Z_{d-i-1} \leq r-j$ on $T$,
and yielding the inequality $a^Y_j + a^Z_{r-j} \geq d-1$ of statement (ii).

Finally, for statement (iii), we need only combine this argument with our
earlier observation that at the closed point, where by hypothesis we had
$a^Y_j+ a^Z_{r-j} = d$ for all $j$, when $i=a^Y_j$ we necessarily have 
$\dim V^Y_i + \dim V^Z_i = r+2$; thus we can in fact choose a basis of 
$V_i$ (still modulo $\m$) which has rank $r+1-j$ in $V^Y_i$ and rank $j+1$ 
(rather than $j$) in $V^Z_i$, and we then get the desired inequality
$a^Y_j + a^Z_{r+1-j} \geq d$ for the entire $T$-valued point, as desired.
Note that the subscheme of $G^r_d$ in question is open simply because
$G^r_d$ is known to map set-theoretically into $G^r_{d,\EH}$, and 
$G^{r, \refn}_{d, \EH}(X)$ is open inside $G^r_{d,\EH}$.
\end{proof}

We observe that for a crude Eisenbud-Harris limit series, there may be many
ways of filling in the intermediate $V_i$ from $V_0$ and $V_0$, so the fiber
of $FR$ may be positive-dimensional. However, the situation is easier to 
get a handle on for the open subset of 
refined series which Eisenbud and Harris actually used in their construction.
Indeed, we show that the space of refined limit series is actually isomorphic 
to an open subscheme of our $G^r_d$ scheme. 

\begin{prop}\label{grd-exact-refined}
Suppose that $(\L^Y, V^Y)$ and $(\L^Z, V^Z)$ form a $T$-valued point of 
$G^{r, \refn}_{d, \EH}(X/B)$. Then we have that
$(\L^Y, V^Y)$ and $(\L^Z, V^Z)$ are the image of a unique $T$-valued point 
under $FR$.
\end{prop}

\begin{proof}It clearly suffices to handle the case that $T$ is connected,
so we make this hypothesis. In this case, we see that we get unique 
vanishing sequences at $\Delta'$ for $V^Y$ and $V^Z$, in the sense that some sequence is satisfied everywhere on $T$, with no stronger ramification index satisfied anywhere on $T$. Indeed, the subscheme of $T$ satisfying $a^Y_i(\Delta') + a^Z_{r-i}(\Delta') > d$ is empty by hypothesis, and because ramification conditions are closed, the ramification sequences obtained at any point of $T$ persist in an open and closed neighborhood.
Now, if $\beta^Y_i$ is the evaluation map $V^Y \rightarrow \pi_{T*}
\L^Y|_{i\Delta'}$ and similarly for $\beta^Z_i$,
this immediately implies that each $\beta^Y_i$ and $\beta^Z_i$ 
has rank determined exactly by the vanishing sequences, in the strong sense
that for some $j$, the closed subscheme where the rank is less than or equal 
to $j$ is all of $T$, but the closed subscheme where the rank is strictly
less than $j$ is empty. It follows (see, e.g., \cite[Prop. 20.8]{ei1}) that the images of the $\beta^Y_i$ and $\beta^Z_i$ are locally free, with
locally free quotients. If we denote by $\{a_j\}_j$ the vanishing sequence 
at $\Delta'$ for $V^Y$, we also note that $\ker \beta^Y_i$ will have
rank $r+1-j$ if $a_{j-1} < i \leq a_j$, and following through the
definitions we see that $\ker \beta^Z_{d-i}$ 
will have rank $j$ if $a_{j-1} \leq i < a_j$, so we find that
$\ker \beta^Y_i = \ker \beta^Y_{i+1}$ if and only if
$\ker \beta^Z_{d-i} = \ker \beta^Z_{d-i+1}$, and 
$\rk \ker \beta^Y_i + \rk \ker \beta^Z_{d-i+1} = r+1$ for all $i$.

The main idea is to construct the $V_i$ as the
subspace of $\ker \beta^Y_i \oplus \ker \beta^Z_{d-i}$ which agree on the
two maps given by evaluation at $\Delta'$. This would then
be unique by Lemmas \ref{grd-random-easy} and \ref{grd-sub},
so we need only show existence. We work locally on the base, so 
that $\L^Y(-i\Delta')|_{\Delta'} \cong \L^Z(-i'\Delta')|_{\Delta'} \cong \O_{\Delta'} \cong \O_T$ for all $i,i'$, and 
fix a
choice of these isomorphisms. As prescribed for gluing together line bundles
defined on components, we define $\L^i$ by the short exact sequence
$$0 \rightarrow \L^i \rightarrow \L^Y (-i \Delta') \oplus \L^Z
((i-d)\Delta') \rightarrow \O_{\Delta'} \rightarrow 0,$$
and pushforward gives us
$$0 \rightarrow \pi_{T*} \L^i \rightarrow \pi_{T*} \L^Y (-i \Delta')
\oplus \pi_{T*} \L^Z ((i-d)\Delta') \rightarrow \O_T.$$
We then define $V_i$ to be the kernel of the induced map, so that:
$$0 \rightarrow V_i \rightarrow \ker \beta^Y_i \oplus \ker \beta^Z_{d-i}
\rightarrow \O_T.$$

We have to show that $V_i$ is a sub-bundle of $\pi_{T*} \L^i$ of the correct
rank. We first observe that the image $\beta^Y_{i+1} (\ker \beta^Y_i)$ has
locally free quotient in $\pi_{T*} \L^Y|_{(i+1)\Delta'}$, and similarly for
$Z$: this image is inside $\im \beta^Y_{i+1}$ by definition, and the
quotient is easily seen to be isomorphic to $\im \beta^Y_i$, via
the map $\pi_{T*} \L^Y|_{(i+1)\Delta'} \rightarrow \pi_{T*} \L^Y|_{i\Delta'}$.
Thus $\beta^Y_{i+1} (\ker \beta^Y_i)$ is a sub-bundle of a sub-bundle, and
must itself be a sub-bundle of $\pi_{T*} \L^Y|_{(i+1)\Delta'}$. Now, we can
factor $\beta^Y_{i+1}$ restricted to $\ker \beta^Y_i$ as
$$\ker \beta^Y_i \rightarrow \pi_{T*} \L^Y(-i\Delta')|_{\Delta'} \hookrightarrow
\pi_{T*} \L^Y|_{(i+1)\Delta'}$$
and we just showed that the cokernel of the composition is locally free; 
since $\L^Y(-i\Delta')|_{\Delta'}$ is a line bundle,
this means the first map must be either zero or surjective, with
surjectivity precisely when $\rk \ker \beta^Y_i = \rk \ker \beta^Y_{i+1}+1$, 
and
the ranks equal otherwise. We obtain the corresponding result for $Z$, and
immediately conclude that $V_i$ is a sub-bundle of $\ker \beta^Y_i \oplus
\ker \beta^Z_{d-i}$, with equality if and only if both $\ker \beta^Y_i =
\ker \beta^Y_{i+1}$ and $\ker \beta^Z_{d-i} = \ker \beta^Z_{d-i+1}$, and
corank one otherwise. Thus, our hypotheses imply that $V_i$ has rank $r+1$. 
The last observation is that $V_i$ being a sub-bundle of $V^Y \oplus V^Z$
implies
that it is a sub-bundle (in our generalized sense) of $\pi_{T*} \L^i$, but this
follows easily from the fact that $V^Y$ and $V^Z$ are sub-bundles of $\pi_{T*}
\L^Y$ and $\pi_{T*} \L^Z$.
\end{proof}

We immediately conclude:

\begin{cor}\label{grd-eh2}The map 
$FR: G^r_d(X/B) \rightarrow G^r_d(Y/B) \times G^r_d(Z/B)$ 
induces an isomorphism from an open subscheme 
$G^r_d(X/B)$ onto $G^{r, \refn}_{d,\EH}(X/B)$, and on the corresponding
separable subschemes of these.
\end{cor}

We may therefore think of the scheme of refined Eisenbud-Harris limit series
as forming an open subscheme of our $G^r_d$ scheme itself:

\begin{defn}We say that a point of $G^r_d$ is a {\bf refined} limit
series if it maps under $FR$ to $G^{r,\refn}_{d,\EH}$, and we denote the 
open subscheme of refined limit series by $G^{r,\refn}_d \subset G^r_d$.
\end{defn}

We also have the following trivial observation.

\begin{cor} Lemma \ref{grd-fr}, Propositions \ref{grd-eh1} and \ref{grd-exact-refined}, and Corollary \ref{grd-eh2} all hold when ramification conditions are imposed.
\end{cor}

\begin{proof}Indeed, we specified ramification solely 
on $V_0$ or $V_d$ depending on whether the relevant section was on $Y$ or
$Z$, so the ramification conditions are visibly compatible with $FR$.
\end{proof}

Since in practice it is less cumbersome to work with Eisenbud-Harris 
series on a given reducible curve, we state our main corollary for the
finite case of Theorem \ref{grd-main} in a situation where one can (nearly)
restrict attention entirely to the Eisenbud-Harris series. We now drop the
hypothesis that we are in case (2), and for notational convenience define:

\begin{defn} In case (1), i.e., with $X/B$ smooth, we simply define 
$G^r_{d, \EH}(X/B)$ to be equal to 
$G^r_d(X/B)$ and similarly for $G^{r,\sep}_{d,\EH}(X/B)$ and $G^{r,\refn}
_{d,\EH}(X/B)$.
\end{defn}

\begin{cor}\label{grd-finite-cor}In the situation of Theorem \ref{grd-main}, 
suppose that $B=\Spec A$ with $A$ a 
DVR having algebraically closed residue field, and $\rho=0$, and denote by
$X_0$ and $X_1$ the special and generic fibers of $X/B$. We omit the
ramification conditions from our notation for the sake of brevity. Then 
consider the following conditions:
\begin{Ilist}
\itm $G^{r,\sep}_{d,\EH}(X_0) \subset G^{r, \refn}_{d, \EH}(X_0)$
\itm $G^{r,\sep}_{d,\EH}(X_0)$ consists of $m$ reduced points for some
$m>0$.
\itm For any DVR $A'$, and any $A'$-valued point of $G^r_d(X)$
such that the induced map $\Spec A' \rightarrow \Spec A$ is flat and the
generic point of $\Spec A'$ maps into $G^{r,\sep}_d(X)$, then the closed
point of $\Spec A'$ maps into $G^{r,\sep}_d(X)$ as well.
\end{Ilist}
If (I) and (II) hold, we have that the $G^{r, \sep}_d(X_1)$ geometrically
contains at least $m$ points; if further (III) holds, then 
$G^{r,\sep}_d(X)$ is finite \'etale over $B$, and the
geometric generic fiber also consists of exactly $m$ reduced points.
\end{cor}

\begin{proof}First, we have by virtue of (I) and Corollary \ref{grd-eh2} 
that $G^{r,\sep}_{d,\EH}(X_0)
\cong G^{r,\sep}_d(X_0)$. Setting $U=G^{r,\sep}_d(X/B)$, if we choose any 
point $x$ in the special fiber, 
applying Corollary \ref{grd-main-cor} we find 
that any component $Z$ of $U$ passing through $x$ maps dominantly to $B$
with $0$-dimensional generic fiber. To count the number of
points, we can take $Z$ to be reduced, in which case it is  
flat over $B$, and we obtain the first assertion.

In the case that (III) holds, we claim that $U$ is in fact
proper over $B$. We begin by noting that the generic fiber must be
$0$-dimensional: by the preceding arguments, this will follow if we show 
that every component $Z$ of $U$ meets the special fiber of $U$, and this
follows from the properness of $G^r_d(X/B)$ and condition (III). 
Now, it suffices to show that $U$ is closed in $G^r_d(X)$, so
choose $y \in U$, $y' \in G^r_d(X)$ distinct points with $y$ specializing to
$y'$; by \cite[Prop. 7.1.9]{ega2} we can find a DVR $A'$ and a map $\Spec
A' \rightarrow G^r_d(X)$ with the generic point mapping to $y$ and the
special point mapping to $y'$. The image of $\Spec A'$ cannot be contained 
in either the special or generic fiber by $0$-dimensionality. Therefore, it 
gives a flat map $\Spec A' \rightarrow \Spec A$, and 
hypothesis (III) implies $y' \in U$ as well, yielding properness of $U$.
Given this, since $U$ is 
unramified in the special fiber, it must be unramified over $B$;
thus, the fibers are reduced, and the lemma which follows gives flatness, 
so we conclude the desired finite \'etaleness.
\end{proof}

The first statement of the following lemma was provided by Max Lieblich.

\begin{lem}\label{aux-reduced} Let $f:X \rightarrow Y$ be a morphism, with
all fibers of $f$ reduced. Then if $f':X_{\red} \rightarrow Y$ is flat, we
have that $X$ is reduced. In particular, if $X$ is irreducible, $Y = \Spec
A$ for some DVR $A$, and $f$ is dominant (still having reduced fibers), then
$f$ is flat.
\end{lem}

\begin{proof} The first assertion follows from Nakayama's lemma together
with consideration of the exact sequence 
$$0 \rightarrow \sN_X \rightarrow \O_X \rightarrow \O_{X_{\red}} \rightarrow
0,$$
where $\sN_X$ is the sheaf of nilpotents inside $\O_X$. For the second
assertion, we apply the standard criterion for flatness over a DVR twice:
first to $X_{\red}$ so we can use the lemma to conclude that $X$ is reduced, 
and then again to $X$ itself.
\end{proof}

\begin{rem}\label{grd-ex1} Note that even if $G^r_d$ maps scheme-theoretically into $G^r_{d
,\EH}$, the statement of Proposition \ref{grd-eh1} would be false if we 
replaced
$d-1$ in the inequality by $d$: we would expect it to fail precisely at the
intersection of the closed subschemes defined by different choices of
$a^Y_j$ and $a^Z_j$ with $a^Y_j +a^Z_{r-j} =d$. Indeed, it is easy enough to write down examples where one does not have the desired inequality:
the simplest case is $Y\cong Z\cong \P^1$, with affine coordinate
functions $y$ and $z$ vanishing at the node, $d=2$, $r=0$, and a
$\Spec k[\epsilon]/(\epsilon^2)$-valued limit series given by 
$$(y^2 + \epsilon y, 0), (y+ \epsilon, \epsilon), (\epsilon, z + \epsilon).$$
Here we are identifying sections of $\O(d')$ with polynomials of degree $d'$,
and for each $i$ specifying pairs of sections of degree $2-i$ and $i$ on the
components, required to agree at $y=z=0$; the inclusion maps are then given
by multiplication by $y$ or $z$ on the appropriate component, and $0$ on the 
other component. This has $a^Y_0 =1$, $a^Z_0=0$; note that
modulo $\epsilon$, it has $a^Y_0 = 2$, $a^Z_0=1$, so it does not correspond
to a refined series, as is required by the proof of Proposition 
\ref{grd-exact-refined} (iii).
\end{rem}

\begin{rem} Continuing along the same line of reasoning, we see that the
scheme-theoretic statement of Proposition \ref{grd-eh1} is actually 
surprisingly
strong; indeed, it implies that if there are two components of the locus of
refined series which meet at a point of $G^r_d$, then if the components have
vanishing sequences at the node given by $a^Y_j, a^Z_j$, and $a'^Y_j,
a'^Z_j$, we must have $|a^Y_j - a'^Y_j| \leq 1, |a^Z_j - a'^Z_j| \leq 1$ for
all $j$.
\end{rem}

\begin{rem} Note that it is not hard to deduce from the proof of
Proposition \ref{grd-eh1} than any refined point is the image of an exact 
point of the linked
Grassmannian used in the construction (see Definition \ref{grd-lg-exact}).
However, the converse is false. Indeed, there exist non-refined points for 
which there is an exact point above them, and there exist others for which
there isn't. We see both already in the simplest case of 
$Y \cong Z \cong \P^1$, and $d=2, r=0$. In the notation of Remark
\ref{grd-ex1}, we could consider $(y^2, 0)$, $(y, z)$, $(0, z^2)$. One
checks that this is not refined, but is an exact point.
On the other hand, if we start with the two 
pairs $(y^2, 0)$, $(0, z^2+z)$, it is easy to see that up to scalar the 
only way to fill in the middle pair is with $(y,0)$, and this point is not 
exact.
\end{rem}

\section{Further Questions}\label{s-grd-further}

This construction brings up a number of natural further questions, and we 
briefly set out a few of them. First, as mentioned earlier, the
Eisenbud-Harris limit series scheme on a reducible curve was never
connected. However, in our case it seems as though the crude limit series
ought to serve as bridges between components of refined limit series with
differing ramification sequences at the node. In fact, at first blush it may 
appear based on dimension-counting that crude limit series should simply be 
the closure of the refined series in many cases, and this may well be true 
in the Eisenbud-Harris scenario of only looking at a $\g^r_d$ on each
component, but because our crude series will often map with 
positive-dimensional fibers to the Eisenbud-Harris crude series, the geometry 
is not entirely clear. For similar reasons, even though our construction
{\it a priori} gives results on smoothing of crude series, the expected 
dimension
hypothesis for all limit series will not follow immediately from having the 
expected dimension of refined series. We can
therefore reasonably ask:

\begin{ques}When is the space of limit series on a reducible curve
connected?
\end{ques}

\begin{ques}When is the space of refined limit series dense in the space of
all limit series?
\end{ques}
 
\begin{ques}In characteristic $0$, what can we say about the dimension of 
spaces of crude limit 
series, and by extension their smoothability? In particular, can we smooth a
``general'' crude series, as we can in the case of refined series (the
latter follows from \cite[Thm.  4.5]{e-h1})?
\end{ques}

We remark that bounding the dimension of crude series on a reducible curve, 
given an understanding of dimensions of linear series on each component,
should be a combinatorial problem, and if the bound is restrictive enough to
imply that on a general curve the crude series have dimension at most as
large as the dimension of refined series, it will follow that for a general
reducible curve, we can always apply the strong form (that is, part (ii))
of Corollary \ref{grd-main-cor} to our entire $G^r_d$ space. In particular,
we can actually make use of the properness of the constructed $G^r_d$ scheme
to obtain direct arguments for theorems such as Brill-Noether, without
requiring arguments involving blowing up the family, as used in \cite[p.
261]{h-m}.

Given our inability to adequately describe the $T$-valued points of $G^r_{d,
\EH}(X/B)$, we can also ask:

\begin{ques}Does $G^r_d(X/B)$ actually map scheme-theoretically into
$G^r_{d, \EH}(X/B)$? Is it scheme-theoretically surjective?
\end{ques}

In applications, an important direction of
generalization is specified ramification along one or more unspecified
smooth sections; this may now be accomplished just as with the case of the 
Eisenbud-Harris theory, by looking at positive-dimensional 
``special fibers'' and allowing $\rho$ to become negative; see
\cite[p. 270]{h-m} for an example.

Finally, the transparency of the construction presented here offers various 
possibilities for generalization beyond the setting of linear series on
curves of compact type. One direction of generalization is to replace 
curves by higher-dimensional varieties. To carry out our main theorem 
in this setting seems at this point to be just a formality, but its
application presents considerable challenges, the most formidable of which
is that the ``expected dimension''
hypothesis of our main theorem is suddenly more of a burden in dimension
higher than one. This is amply demonstrated by the interpolation problem
(see \cite{ci1} and \cite{gi3}), where one sees first that expected 
dimension for general ramification points need not hold, even for 
zero-dimensional linear series on $\P^2$, and second, that standard 
degeneration arguments have thus far failed to succeed in describing when 
exactly the expected dimension is in fact correct.

One could also hope to generalize from line bundles to vector bundles. This
should not be too difficult, but in order to set up inductive degeneration
arguments, one would then need some a description of the limit objects on 
the reducible curve
which would play the role of the Eisenbud-Harris description of limit 
series. Lastly, it might be possible to adapt the construction to work on
curves not of compact type, but in this setting one may find that the
functor for a given family, after restriction to the reducible special
fiber, will still depend on the geometry of the entire family. This would 
potentially complicate the situation considerably.

\appendix 

\section{The Linked Grassmannian Scheme}\label{s-grd-lg}

In this appendix, we develop of a theory of a moduli scheme parametrizing
collections of sub-bundles of vector bundles on a base scheme, linked
together via maps between the vector bundles. Representability by a proper
scheme is easy and true quite generally; however, to obtain dimension
formulas will require more hypotheses and more work. These
hypotheses, while reasonably natural and easy to state, are motivated by the
idea that the vector bundle maps are induced as pushforwards of certain 
maps of sufficiently ample line bundles on a scheme proper over the base 
scheme, as in the situation of the limit linear series theory of the present
paper, and its natural generalization to higher-dimensional varieties. 

We first specify the objects we will study in more detail; for the remainder
of this appendix, we will be in: 

\begin{sit}\label{grd-lg-sit}Let $S$ be any base scheme, and 
$\E_1, \dots, \E_n$ vector
bundles on $S$, each of rank $d$. We have maps 
$f_i:\E_i \rightarrow \E_{i+1}$ and $g_i: \E_{i+1} \rightarrow
\E_i$, and a positive integer $r<d$.
\end{sit}

The functor we wish to study may now be easily described:

\begin{defn}\label{grd-lg-functor}In this situation, we have the functor 
$\cLG(r,\{\E_i\}_i,
\{f_i, g_i\}_i)$, associating to each $S$-scheme $T$ the set of sub-bundles
$V_1, \dots, V_n$ of $\E_{1,T}, \dots, \E_{n,T}$ of rank $r$ and satisfying 
$f_{i,T}(V_i) \subset V_{i+1}$, $g_{i,T}(V_{i+1}) \subset V_i$ for all $i$.
\end{defn}

Then without any further hypotheses, we have:

\begin{lem}\label{grd-pre-lg}$\cLG(r,\{\E_i\}_i, \{f_i, g_i\}_i)$ is representable by 
a projective scheme $\LG$ over $S$, which is naturally a closed subscheme of
a product $G$ of Grassmannian schemes over $S$; $G$ is smooth and 
projective over $S$ of relative dimension $n r(d-r)$.
\end{lem}

\begin{proof}
Let $G_i$ be the schemes of Grassmannians of rank $r$ sub-bundles of the
$\E_i$, and $G$ the product of the $G_i$ over $S$. Denote our projection 
maps from $G$ to each $G_i$ by $\pi_i$, and the maps from each $G_i$
to $S$ by $\phi_i$. Let $\F_i$ be the universal sub-bundles on each $G_i$. 

Then each $f_i$ induces a map 
$$\pi_i ^* \F_i \rightarrow \pi_i^* \phi_i^* \E_i = \pi_{i+1} ^* \phi_{i+1}^*
\E_i \stackrel{f_i}{\rightarrow} \pi_{i+1} ^* \phi_{i+1}^* \E_{i+1} 
\rightarrow \pi_{i+1} ^* \phi_{i+1} ^* \E_{i+1} / \pi_{i+1}^* \F_{i+1}$$
on $G$, and 
the kernel of this map is a closed subscheme which imposes precisely the
condition that $f_i(V_i) \subset V_{i+1}$. Similarly, $g_i$ induces a map
$\pi_{i+1} ^* \F_{i+1} \rightarrow \pi_i ^* \phi_i ^* \E_i / \pi_i^* \F_i$ on 
$G$ whose kernel imposes the condition $g_i(V_{i+1}) \subset V_i$. Taking
the intersection of these closed subschemes for all $f_i$ and $g_i$ thus
gives a scheme representing $\cLG(r,\{E_i\}_i, \{f_i, g_i\}_i)$, which as a
closed subscheme of $G$ is projective over $S$.
\end{proof}

However, in order to say anything of substance about the $LG$ scheme itself, 
and in particular to get the necessary lower bound on
dimension, we need to make a number of additional hypotheses. We define:

\begin{defn}\label{grd-lg-def}In Situation \ref{grd-lg-sit}, we say that 
$\LG(r,\{E_i\}_i, \{f_i, g_i\}_i)$ is a {\bf linked 
Grassmannian} of {\bf length $n$} if $S$ is integral and Cohen-Macaulay, and 
the following 
additional conditions on the $f_i$ and $g_i$ are satisfied:

\begin{Ilist}
\itm There exists some $s \in \O_S$ such that $f_i g_i = g_i f_i$ is
scalar multiplication by $s$ for all $i$.
\itm Wherever $s$ vanishes, the kernel of $f_i$ is precisely equal to
the image of $g_i$, and vice versa. More precisely, for any $i$ and given any 
two integers $r_1$ and $r_2$ such that $r_1 + r_2 < d$, then the closed
subscheme of $S$ obtained as the locus where $f_i$ has rank less than or
equal to $r_1$ and $g_i$ has rank less than or equal to $r_2$ is empty.
\itm At any point of $S$, $\im f_i \cap \ker f_{i+1}=0$, and $\im
g_{i+1} \cap \ker g_i = 0$. More precisely, for any integer $r_1$, and any
$i$, we have locally closed subschemes of $S$ corresponding to the locus
where $f_i$ has rank exactly $r_1$, and $f_{i+1} f_i$ has rank less than or
equal to $r_1-1$, and similarly for the $g_i$. Then we require simply that 
all of these subschemes be empty.
\end{Ilist}
\end{defn}

\begin{rem}The hypothesis that $S$ is integral and Cohen-Macaulay is 
unnecessary for most of our analysis. We use it only in the 
dimension theory portion of the argument, to ensure that $\LG$ is catenary.
\end{rem}

From this point on, we strengthen Situation \ref{grd-lg-sit}.

\begin{sit} We suppose that $\LG$ is a linked Grassmannian, and we
denote its structure map to $S$ by $\pi$.
\end{sit}

The following lemma will be convenient for constructing and manipulating
points of $\LG$:

\begin{lem}\label{grd-lg-decomp}Let $\{V_i\}_i$ be a k-valued point of $\LG$, and 
suppose $s=0$ in
$k$. Then for any $i$, we can decompose $V_i$ as $f_{i-1}(V_{i-1}) \oplus
\ker f_i |_{V_i} \oplus C$ for some complementary subspace $C \subset V_i$.
Indeed, if we specify any $C' \subset \ker g_{i-1}|_{V_i}$ which intersects
$f_{i-1}(V_{i-1})$ trivially, we may choose $C = C' \oplus C''$ for some
$C''$.
\end{lem}

\begin{proof}Clearly, it suffices to show that for any $C'$ as in the
statement, we have that $f_{i-1}(V_{i-1})\oplus \ker f_i|_{V_i} \oplus C'$ 
injects into $V_i$. But suppose we have $v_1 \in f_{i-1}(V_{i-1})$, $v_2 \in
\ker f_i|_{V_i}$, and $v_3 \in C'$, such that $v_1+v_2+v_3 =0$. If we apply
$g_{i-1}$, by hypothesis $g_{i-1}(v_3)=0$, and $g_{i-1}(v_1)=0$ because $v_1$ 
is in the image of $f_{i-1}$ and we assumed $s=0$. So we find that 
$g_{i-1}(v_2)=0$, which we claim implies $v_2=0$: indeed, $v_2 \in \ker f_i$
by hypothesis, so by condition (II) of a linked Grassmannian it is in
the image of $g_i$, and by condition (III), it cannot map to $0$ under
$g_{i-1}$ unless it is $0$. Hence $v_2=0$, so $v_1+v_3=0$, and since we
assumed that $C'$ was disjoint from $f_{i-1}(V_{i-1})$, we get $v_1=v_3=0$
as well.
\end{proof}

In order to make inductive arguments convenient, we define:

\begin{defn}If $\LG$ is a linked Grassmannian of length $n$, and $n'$ any
positive integer less than $n$, we have the {\bf truncation map} from
$\LG$ to the linked Grassmannian of length $n'$ obtained by forgetting all
$\E_i$, $f_i$, and $g_i$ for all $i>n'$.
\end{defn}

We will want to know that the truncation map is always surjective, even on
certain classes of families:

\begin{lem}\label{grd-lg-lift}The truncation map is surjective for all $n'$. 
Further, in the case that the base is a point, let $x = \{V_i\}_i$ be any
point of $\LG$, and suppose we have a family 
$\tilde{x}_{n'} = \tilde{V}_i|_{i \leq n'}$ (that is, a scheme-valued point 
of the restricted $LG$ scheme) specializing to the truncation of 
$x$ to length $n'$, and such that $\tilde{V}_{n'}$ may be written as 
$\tilde{C}_{n'} \oplus \ker f_{n'}|_{V_{n'}}$ for some family 
$\tilde{C}_{n'}$. Then
$\tilde{x}_{n'}$ may be lifted to a family $\tilde{x}$ of length $n$, 
specializing to $x$, possibly after a Zariski localization of the base of
the family.
\end{lem}

\begin{proof}Surjectivity may be checked on points, and given the description 
of the $\cLG$ functor, it suffices to handle the case $n=2$, $n'=1$.
Over a point, we may consider $\E_1 =
\E_2 = E$ to be a single $d$-dimensional vector space, and $f_1$ and $g_1$ to 
be self-maps of $E$. Let $V_1$ be a vector space of dimension $r$ inside $E$;
we just need to show that there exists a $V_2$ of dimension $r$ inside $E$
such that $f_1(V_1) \subset V_2$, and $g_1(V_2) \subset V_1$, or
equivalently, such that $f_1(V_1) \subset V_2 \subset g_1^{-1}(V_1)$.
Now, $\dim f_1(V_1) \leq r$ and $f_1(V_1) \subset 
g_1^{-1}(V_1)$ by hypothesis, so it suffices to observe
that $\dim g_1^{-1}(V_1) = \dim \ker g_1 + \dim (V_1 \cap \im g_1)$,
and the codimension of $\im g_1$ in $E$ and hence $V_1$ is bounded by 
$\dim \ker g_1$, so we conclude that $\dim g_1^{-1}(V_1) \geq r$. 

For the second assertion, it suffices to show that we can lift to a
$\tilde{V}_{n'+1}$ of the form $\tilde{C}_{n'+1}\oplus \ker
f_{n'+1}|_{V_{n'+1}}$ and specializing to the truncation of $x$ to 
length $n'+1$, since then we can iterate until we have lifted all the way to
length $n$. Thanks to Lemma \ref{grd-lg-decomp}, we can write $V_{n'} =
C_{n'} \oplus \ker f_{n'}|_{V_{n'}}$, and $V_{n'+1}=
f_{n'}(V_{n'}) \oplus \ker f_{n'+1}|_{V_{n'+1}} \oplus C_{n'+1}$ for some
$C_{n'}$ and $C_{n'+1}$, with $\tilde{C}_{n'}$ specializing to $C_{n'}$, and
in particular, having full rank under $f_{n'}$ except possibly on a
closed subset of the base supported away from $x$, where the rank could drop. 
Away from this locus on the base, if we replace $V_{n'+1}$ by
$f_{n'}\tilde{V}_{n'} \oplus \ker f_{n'+1}|_{V_{n'+1}} \oplus C_{n'+1}$
(that is, if we set $\tilde{C}_{n+1} = f_{n'}\tilde{V}_{n'} \oplus
C_{n'+1}$), noting that $f_{n'}\tilde{V}_{n'} = f_{n'}\tilde{C}_{n'}$, we 
obtain a lifting with the desired properties.
\end{proof}

The key notion for getting a handle on the $\LG$ scheme is the following:

\begin{defn}\label{grd-lg-exact}We say that a point of a linked Grassmannian 
scheme is {\bf exact} if the corresponding collection of vector spaces $V_i$
satisfy the conditions that $\ker g_i |_{V_{i+1}}\subset f_i(V_i)$ and 
$\ker f_i|_{V_i} \subset g_i(V_{i+1})$ for all $i$.
\end{defn}

The last part of assertion (ii) of the following lemma is gratuitous, but it
follows immediately from the argument for the rest, and may perhaps
shed some little light on the overall situation.

\begin{lem}\label{grd-lg-exact-desc}We have the following description of 
exact points:
\begin{ilist}
\itm The exact points form an open subscheme of 
$\LG$, and are naturally described as the complement of the closed subscheme 
on which $\rk f_i|_{V_i}+ \rk g_i|_{V_{i+1}} < r$ for some $i$. 
\itm In the case $s=0$, we find that we can describe exact points as 
those with
$\rk f_i(V_i)+ \rk g_i(V_{i+1}) = r$ for all $i$, even for arbitrary 
scheme-valued points, and we also find that an exact point has $f_i(V_i)$ a
sub-bundle of $V_{i+1}$, and $g_i(V_{i+1})$ a sub-bundle of $V_i$ for all 
$i$.
\end{ilist}
\end{lem}

\begin{proof}We certainly get a closed subscheme as described, simply by
taking the union over all $i$ and $r_1, r_2$ with $r_1+r_2<r$ of the loci
described by $\rk f_i|_{V_i} \leq r_1$ and $\rk g_i|_{V_{i+1}} \leq r_2$. We
immediately see that the points of this set are precisely the complement of
the exact points, since outside the locus where $s$ vanishes, both $f_i$ and 
$g_i$ are invertible, and correspondingly all points are exact; on the other 
hand, if
$s$ vanishes at our point, we have $\im f_i \subset \ker g_i$ and 
$\im g_i \subset \ker f_i$ for all $i$, so we already have that
$\dim f_i(V_i) + \dim g_i(V_{i+1}) \leq r$ and we get strict inequality if
and only if these containments are strict. 

For the second part, if a $T$-valued point 
satisfies $\rk f_i|_{V_i}+ \rk g_i|_{V_{i+1}} = r$ for all $i$ on all of $T$,
by definition it is in the complement of the closed subscheme of non-exact
points defined above, so it is certainly exact. Conversely, for the other
direction it suffices to work over local rings, so suppose we
have a $T$-valued point $(V_i)$ of $\LG$ where $T$ is local, and the point
$(\bar{V}_i)$ of $\LG$ at the closed point of $T$ is exact. 
Since $s$ is zero on $T$, then for any given $i$ we have $\ker \bar{f}_i
|_{\bar{V}_i}= \bar{g}_i(\bar{V}_{i+1})$ and vice versa; in particular, if 
we choose 
$\bar{v}_1, \dots, \bar{v}_{r_1}\in \bar{V}_i$ such that the 
$\bar{f}_i(\bar{v}_j)$ form a basis of $\bar{f}_i(\bar{V}_i)$, and 
$\bar{v}'_1, \dots, \bar{v}'_{r_2}\in\bar{V}_{i+1}$ such that the
$\bar{g}_i(\bar{v}'_j)$ form a basis of $\bar{g}_i(\bar{V}_{i+1})$, we 
find that $r_1+r_2=r$, and we obtain a basis $\bar{e}_i$ (resp., 
$\bar{e}'_i$) for $\bar{V}_i$ (resp., $\bar{V}_{i+1}$) given by 
the $\bar{v}_j$ and $g_i(\bar{v}'_j)$ (resp., $\bar{v}'_j$ and 
$f_i(\bar{v}_j)$). By Nakayama's lemma, we can lift this situation to the
local ring, and easily check that the desired assertions follow.
\end{proof}

Our main technical lemma for this appendix is:

\begin{lem}\label{grd-lg-exact-props} We have the following statements on exact
points: 

\begin{ilist}
\itm The exact points are dense in $\LG$, and indeed dense in every
fiber.
\itm Given any exact point $x \in \LG$, let $y$ be its image in $S$,
suppose $A$ is a local ring, and $A'$ a quotient of $A$. Let $T= \Spec A$,
and $T' =\Spec A'$. Then given any commutative diagram containing the solid 
arrows of
$$\xymatrix{{T'} \ar^-{f}[r] \ar[d] & {\LG} \ar[d] \\ {T} \ar[r] \ar@{-->}[ur] &
{S}}$$
with the closed point of $T'$ mapping to $x$, the dashed arrow may also be 
filled in. In particular, $x$ is a smooth point of $\LG$ over $S$.
\end{ilist}
\end{lem}

\begin{proof} For (i), 
To see that the exact points are dense in every fiber, suppose we have a 
non-exact
point; we just observed that this corresponds to a set of $V_i$ such that for 
at least one $i$, we have $\dim f_i(V_i)+ \dim g_i(V_{i+1}) < r$.
In particular, we are in the situation where $f_i g_i = g_i f_i =0$. 
Now, choose the 
smallest $i$ such that $\dim f_i(V_i)+ \dim g_i(V_{i+1}) < r$, and truncate 
our linked Grassmannian to $i+1$; 
here, we show that there are nearby points in the fiber such that
the condition $\dim f_i(V_i)+ \dim g_i(V_{i+1}) = r$ is satisfied. We leave
$V_1$ through $V_i$ unmodified. By hypothesis, there are vectors in
$V_{i+1}$ in the kernel of $g_i$ which are not in $f_i(V_i)$, and
vice versa; indeed, we see that $r':= \dim \ker g_i|_{V_{i+1}}-\dim f_i(V_i) =
\dim \ker f_i|_{V_i}-\dim g_i(V_{i+1})=r - \dim f_i(V_i)-\dim
g_i(V_{i+1})$.
Choose $C_i$ and $C_{i+1}$ in $\ker f_i|_{V_i}$ and $\ker g_i|_{V_{i+1}}$ of
dimension $r'$, intersecting $g_i(V_{i+1})$ and $f_i(V_i)$ trivially; we
have that together with these spaces, they must complete
the span of $\ker f_i|_{V_i}$ and $\ker g_i|_{V_{i+1}}$ respectively. 
Since $C_i \subset \ker f_i|_{V_i}$, it is in $\im g_i$, and we can find
$e_1, \dots, e_{r'} \in \E_{i+1}$, whose span is necessarily disjoint from 
$V_{i+1}$, and which map to a 
basis of $C_i$ under $g_i$. By Lemma \ref{grd-lg-decomp}, we can write 
$V_{i+1}=f_i(V_i)\oplus \ker
f_{i+1}|_{V_{i+1}} \oplus C_{i+1} \oplus C''$ for some $C''$. If we take any
basis $e'_1, \dots, e'_{r'}$ for $C_{i+1}$, we can make a family
$\tilde{V}_{i+1}$ over $\A^1$ by replacing $C_{i+1}$ with the span of $e'_i
+ t e_i$ for all $i$, as $t$ varies. 

Now, $\tilde{V}_{i+1}$ specializes to
$V_{i+1}$ at $t=0$, and we see that it always remains linked to $V_1, \dots,
V_i$, left unmodified: it certainly maps into $V_i$ under $g_i$, since we
are modifying basis elements by the $e_i$, which were chosen to map into
$V_i$; on the other hand, our construction leaves the summand
$f_i(V_i)$ unmodified, so $f_i$ certainly maps $V_i$ into any member of
$\tilde{V}_{i+1}$. We also observe that we now have $\dim
f_i(V_i)+ \dim g_i(\tilde{V}_{i+1})=r$ whenever $t \neq 0$: indeed, 
$C_{i+1}$ was in the kernel of $g_i$ for $t=0$, so we still have
$g_i(\tilde{V}_{i+1}) \supset g_i(V_{i+1})$; for any
$t \neq 0$, $C_{i+1}$ maps isomorphically to $C_i$ under $g_i$; finally,
since we chose $C_i$ to, together with $g_i(V_{i+1})$, span 
$\ker f_i|_{V_i}$, we find that for any $t \neq 0$, $g_i(\tilde{V}_{i+1}) =
\ker f_i|_{V_i}$, giving the desired exactness at $i$. Now, by Lemma
\ref{grd-lg-lift}, we can lift this family to a family $\tilde{V}_j$ for all
$j$, specializing to our given point, but now satisfying 
$\dim f_i(V_i)+ \dim g_i(V_{i+1}) = r$ for a general point in the family; we
conclude that the points which are non-exact at the $i$th step (but exact
for $j<i$) are in the closure of those which are exact through the $i$th
step, and by induction are actually in the closure of the points which are
exact at all steps.

For assertion (ii), $f(T')$ corresponds to a collection $\{V_i\}_i$ over $A'$;
Our $\E_i$ are now all free modules of rank $d$ over $A$, and we simply want
to produce free $A$-submodules $\tilde{V}_i$ (with free quotients) linked 
by the $f_i$ and $g_i$
and restricting to the given $V_i$ in the quotient ring $A'$. 
To do this, denote by $\bar{V}_i$ the collection of subspaces over
the residue field of $A'$ corresponding to $x$, and let $r_i, r_i'$ be the 
dimensions of $f_i(\bar{V}_i), g_i(\bar{V}_{i+1})$ respectively for each $i$. 
We begin by choosing bases $\bar{e}^i_j$ of $\bar{V}_i$, and lifting
appropriately. If our $f_i$ and $g_i$ are invertible at the closed point, 
which is to say, if the $s$ from 
condition (I) of a linked Grassmannian is non-zero in $\kappa(y)$, we simply 
choose an arbitrary basis $\bar{e}^1_j$ of $\bar{V}_1$, and take its images 
under the $f_i$. We then lift the $\bar{e}^1_j$ to $V_1$, and take images
under the $f_i$, to obtain bases of the $V_i$, and lift by the same process
to $\E_i$, defining sub-modules $\tilde{V}_i$.

Otherwise, if we had $s=0$ in $\kappa(y)$, for each $i$ we choose 
$\bar{e}^i_j$ in three 
categories: first, $r-r_{i-1}-r'_i$ elements which are linearly independent
from the span of $f_{i-1}(\bar{V}_{i-1})\cup g_i(\bar{V}_{i+1})$; second,
$r_{i-1}$ elements generating $f_{i-1}(\bar{V}_{i-1})$; and third,
$r'_i$ elements generating $g_i(\bar{V}_{i+1})$. Noting that 
even without exactness, since $s=0$, we have $r = \dim f_{i-1}(\bar{V}_{i-1}) 
+ \dim \ker f_{i-1}|_{\bar{V}_{i-1}} \geq r_{i-1}+r'_{i-1}\geq r_{i-1}+r'_i$, 
and $f_{i-1}(\bar{V}_{i-1}) \subset \ker g_{i-1}|_{\bar{V}_i}$ which is
disjoint from $g_i(\bar{V}_{i-1})$,
so we see that this is possible. Moreover, by choosing the first category 
for all $i$ first, we can inductively construct the basis elements in the 
second and third categories to be images under $f_{i-1}$ and $g_i$ of basis 
elements already chosen, moving from $i=1$ to $i=n$ for the second category, 
and the opposite direction for the third.
Next, choose lifts $e^i_j$ to the $V_i$, using the same process of lifting
all $\bar{e}^i_j$ in the first category first, and defining the rest as
iterated images under $f_{i-1}$ and $g_i$. Finally, lift the $e^i_j$ to 
$\tilde{e}^i_j \in \E_i$, once again via the same process, and define
$\tilde{V}_i$ to be the span of the $\tilde{e}^i_j$. 

By Nakayama's lemma, the $e^i_j$ constructed in either case give free 
generators for the $V_i$. Applying Nakayama's lemma again, we find that 
the $\tilde{V}_i$ are sub-bundles of $\E_i$ of rank $r$, and clearly they 
specialize to the $V_i$, so we need only check that they are linked. In the
case that $s$ was non-zero in $\kappa(y)$, the $\tilde{V}_i$ are linked under 
the $f_i$ by construction, and must likewise be linked under the $g_i$,
since $g_i$ is a unit times the inverse of $f_i$. In the case where $s$ was
zero in $\kappa(y)$, take any $\tilde{e}^i_j$ for $i<n$; we show that its
image under $f_i$ is a scalar multiple of $\tilde{e}^{i+1}_{j'}$ for some
$j'$. We now apply exactness, to note that either $f_i(\bar{e}^i_j)=0$ and 
$\tilde{e}^i_j = g_i(\tilde{e}^{i+1}_{j'})$, or we defined 
$\tilde{e}^{i+1}_{j'}=f_i (\tilde{e}^i_j)$, for some $j'$. In the latter
case, we are done, while in the former case we simply observe that
$f_i (\tilde{e}^i_j) = s (\tilde{e}^{i+1}_{j'})$. The same argument
works for the $g_i$, so we have constructed a map from $T$ to $LG$ lifting 
$f$, which by \cite[Prop.  17.14.2]{ega44} completes the proof of part (ii). 
\end{proof}

The following proposition provides a strong converse to part (ii) of the
above lemma:

\begin{prop}\label{grd-lg-exact-comps}The non-exact points of a fiber are 
precisely the intersections of the components of that fiber.
\end{prop}

\begin{proof}Since the exact points are smooth, they are certainly not in
any intersection of components. For the other direction, we first make the
following observation: because ranks can only drop under specialization, 
given two exact points $\{V_i\}_i$ and $\{V'_i\}_i$, with $r_i:=\dim
f_i(V_i)$ and $r'_i:= \dim f_i(V'_i)$, if some $r_i \neq r'_i$, then the 
two points must lie on distinct components of $\LG$. Thus, to show that any
non-exact point is in the intersection of components, it suffices to exhibit
it as the specialization of two different exact points with distinct $r_i$. 

Looking at the proof of Lemma \ref{grd-lg-exact-props} part (i), we see that any
point which is non-exact at $i_0$, with $i_0$ minimal, can expressed as the 
specialization of an exact point with $r_i$ unchanged for all $i \leq i_0$;
however, upon closer examination, we see that in fact the process leaves all the
$r_i$ unchanged, simply increasing the dimensions of the $g_i(V_{i+1})$ as 
necessary to make the points exact. On the other hand, we note that the
linked Grassmannian situation is completely symmetric in the $f_i$ and 
$g_i$, so now that we have shown that any point can be written as the
specialization of an exact point with the $\dim f_i(V_i)$ unchanged, it
follows by symmetry  that there is another exact point specializing to our 
given point, leaving the dimensions of the $g_i(V_{i+1})$ intact,
and therefore necessarily increasing at least some of the $r_i$. This then
expresses our non-exact point as lying in the intersection of two
components, as desired.
\end{proof}

We can also use the smoothness at exact points to compute the dimension of
fibers of $\LG$:

\begin{lem}\label{grd-lg-fiber-dim}The fibers of $\LG$ over $S$ have every 
component of dimension precisely $r(d-r)$.
\end{lem}

\begin{proof}In view of Lemma \ref{grd-lg-exact-props}, 
we can compute the dimension of any component of the fiber
by showing that its tangent space at any exact point has the desired
dimension. Since we are only looking at a fiber, we set $S=\Spec k$. If $s
\neq 0$ in $k$, $\LG \cong \G(r,d)$, and is smooth of dimension $r(d-r)$, so
there is nothing to show. Otherwise, suppose we have a collection of $V_i$ 
corresponding to an exact point. Then $\ker f_i|_{V_i}=g_i(V_{i+1})$ for all
$i$, so we use Lemma \ref{grd-lg-decomp} to write each $V_i$ as
$f_{i-1}(V_{i-1}) \oplus g_i(V_{i+1}) \oplus C_i$ for some complementary
space $C_i$. Our first assertion is that the dimensions $d_i$ of the $C_i$ 
add up to $r$. Indeed, if we let $r_i = \dim f_i(V_i)$, and 
$r'_i = \dim g_i(V_{i+1})$, we have $r_i = r-r'_i$ from exactness, and for
$1 < i < n$, $d_i = r - r_{i-1}-r'_i = r_i -r_{i-1}$, with 
$d_1 = r-r'_1 = r_1$ and $d_n = r-r_{n-1}$, so we see we indeed have $\sum _i
d_i = r$.

The next claim is that first-order deformations of the $V_i$ inside of $\LG$
correspond precisely to first-order deformations of each $C_i$ individually
inside $\E_i$, taken modulo deformations of the $C_i$ which remain inside 
$V_i$. Any deformation of the $C_i$ together will yield a deformation of 
the $V_i$: we use our direct sum decomposition to inductively define the
induced deformation, obtaining deformations of $f_i(V_i)$ as the image of
the deformation of $C_{i-1}$ together with the (inductively obtained)
deformation of $f_{i-1}(V_{i-1})$, and similarly for the $g_i(V_{i+1})$. 
Moreover, since each $f_i(V_i)$ is spanned by $f_{i-1}(V_{i-1})$ together
with $C_{i-1}$, this is the only possible way to obtain a deformation of the
$V_i$ given deformations of the $C_i$.
Clearly, two deformations of the $C_i$ will yield equivalent deformations of
$V_i$ if and only if their difference is a deformation of the $C_i$ inside 
of its $V_i$. Finally, any deformation of the $V_i$ may be expressed
(non-uniquely) as a deformation of its summands, and in particular gives
a deformation of the $C_i$, at least up to the same equivalence relation.
Since the deformation of the $V_i$ induced by the deformations of the $C_i$
was unique, this must invert our first construction, completing the proof 
of the claim.

Now we are done: first-order deformations of any given $C_i$ are given by
the tangent space to $\G(d_i,d)$, which is a variety smooth of dimension
$d_i(d-d_i)$, so has $d_i(d-d_i)$-dimensional tangent space at any point.
Similarly, the
space of deformations of $C_i$ inside of $V_i$ has dimension $d_i(r-d_i)$; 
the
difference is $d_i(d-r)$. Thus, the total dimension of our tangent space is 
$\sum _i d_i(d-r) = r(d-r)$, as asserted.
\end{proof}

We now have all the tools to prove our main result:

\begin{thm}\label{grd-main-lg}A linked Grassmannian scheme is a closed subscheme
of the obvious product of Grassmannian schemes over $S$; it is
projective over $S$, and each component has codimension $(n-1)r(d-r)$
inside the product, and maps surjectively to $S$. If $s$ 
is non-zero, then $\LG$ is also irreducible.
\end{thm}

\begin{proof}We already have that the linked Grassmannian is projective over
$S$, and lies inside the obvious product of Grassmannians, which we denote
by $G$. It is easy to see 
each component maps dominantly onto $S$, since the exact points are both
smooth and dense by Lemma \ref{grd-lg-exact-props}. 

For the dimension statement,
given any component of $\LG$, let $x$ be an exact point of $\LG$ on the
specified component, and not on any other component, and $s$ the image of $x$
in $S$. Since $S$ is Cohen-Macaulay, everything is catenary, so codimensions
can be computed naively for irreducible spaces. By Lemma 
\ref{grd-lg-fiber-dim}, we have that 
$\O_{\LG,x}$ is smooth over $\O_{S,s}$ of relative dimension $r(d-r)$, and 
in particular integral. Similarly, $\O_{G,x}$ is locally affine over 
$\O_{S,s}$, hence integral and smooth of relative dimension $nr(d-r)$. 
The desired codimension statement then follows from 
\cite[Prop. 17.5.8 (i)]{ega44}.

Finally, when $s$ is non-zero, over the open subset of $S$ where $s$ is
invertible, the fibers are all simply Grassmannians of dimension $r(d-r)$;
since the map is proper, we conclude that $\LG$ is irreducible over
this locus, of dimension $r(d-r)$. On the other hand, since every component
maps dominantly to $S$, there cannot be any component of $\LG$ contained in 
the locus where $s$ vanishes, yielding the desired irreducibility.
\end{proof}

\begin{warn}Lemma \ref{grd-lg-exact-desc} sounds quite innocuous, but there 
are some
pitfalls to be aware of. Consider the simple example of $n=d=2$, $r=1$, 
$S = \Spec k$, $\E_1=\E_2 = k^2$, 
$f_1 = \begin{bmatrix}1 & 0 \\ 0 & 0\end{bmatrix}$, 
and $f_2 = \begin{bmatrix}0 & 0 \\ 0 & 1\end{bmatrix}$. 
In this case, if $V_1$ is generated by 
$v_1= \begin{bmatrix}X_0 \\ X_1\end{bmatrix}$ and $V_2$ by
$v_2= \begin{bmatrix}Y_0 \\ Y_1\end{bmatrix}$, we find the condition for them 
to be linked is simply that $X_0 Y_1 = 0$, and it is easy enough to check that
we actually get that $\LG$ is scheme-theoretically cut out by this equation
inside $\P^1 \times \P^1$, giving a pair of $\P^1$'s attached at $X_0 = Y_1
=0$, which is the only non-exact point. Our lemma has shown that deformations 
have to behave well at the exact points, but if we consider the $T$-valued
point for $T =\Spec k[\epsilon]/(\epsilon^2)$ with $V_1$ generated by
$v_1=\begin{bmatrix}\epsilon \\ 1\end{bmatrix}$ and $V_2$ generated by
$v_2= \begin{bmatrix}1 \\ \epsilon \end{bmatrix}$, we note two pathologies:

First, this point actually satisfies our initial set-theoretic description
of an exact point, that $\ker g_1 |_{V_2}\subset f_1(V_1)$ and vice versa, as 
both images and kernels will be given precisely by $\epsilon v_i$. So
this description, while dealing with both $s=0$ and $s$ invertible
simultaneously, is only valid from a set-theoretic point of view. 

Second, while we have shown that at any (scheme-valued) exact point, there 
will be an $r_1$ and $r_2$ with $r_1 + r_2 = r$ and $\rk f_1|_{V_1} \leq r_1$, 
$\rk g_1|_{V_2} \leq r_2$, we see that by allowing the ranks to drop at the
closed point, we actually allow them to increase on the local ring level.
Specifically, in our case $r=1$, so either $r_1$ or $r_2$ would have to be
$0$, but neither $f_1$ nor $g_1$ is the zero map. Of course, this makes
perfect sense geometrically, as the node will necessarily have tangent 
vectors which don't point along either branch, but it underscores the fact
that the $T$-valued points of a union of schemes is not simply the union of
the $T$-valued points of the individual schemes.
\end{warn}

We conclude with an example and some further questions which we have not
pursued here because they are not necessary for our applications.

\begin{ex}We consider the situation of $S = \Spec k$, $n=2$. In this case, 
it is easy to describe the components explicitly, as well as to see their 
dimensions without invoking any deformation theory. We already know that if
$s \neq 0$, we just get a Grassmannian, so we assume that $s=0$. If we write
$d_1 = \rk f_1$, $d_2 = \rk g_1$ (on the entire vector space), we have $d_1
+ d_2 = d$ by condition (II) of a linked Grassmannian. We will see that
there are $\min\{r+1, d-r+1, d_1+1, d_2+1\}$ components, each of dimension
$r(d-r)$, and indexed by the dimension of $f_1(V_1)$ on general points.

Indeed, we saw in the proof of Lemma \ref{grd-lg-lift} that  
the fiber of any point $V_1$ of $G_1$ under truncation is simply the 
Grassmannian of vector spaces $V_2$ containing $f_1(V_1)$ and contained in
$g_1^{-1}(V_1)$, which had dimension 
$\dim \ker g_1 + \dim (V_1 \cap \im g_1)$. 
We need to see that this dimension depends only on the dimension of 
$f_1(V_1)$, which we
will denote by $r_1$. By condition (II) of a linked
Grassmannian, $\ker g_1 = \im f_1$, and $\im g_1 = \ker f_1$, so we may
write this as $d_1 + \dim (V_1 \cap \ker f_1)$. Furthermore, $\dim (V_1
\cap \ker f_1) = r- r_1$, so we can write everything in terms of $r_1$, as
desired. Specifically, we have a Grassmannian of $r$-dimensional subspaces of 
a $(d_1 + r - r_1)$-dimensional space, containing an $r_1$-dimensional space,
and this has dimension $(r-r_1)(d_1-r_1)$. 

We now obtain our assertions without trouble: fix an 
$r_1 \leq \min\{r,d_1\}$ also satisfying $r_1 \geq \max\{0, r-d_2\}$, and
consider the locally closed subset $G^{r_1}_1$ in $G_1$ with 
$\dim f_1(V_1) = r_1$. Note that the specified range is precisely the range
for which this will be non-empty. Now, $G^{r_1}_1$ is an open subset of the 
locus in $G_1$ with $\dim f_1(V_1) \leq r_1$, which corresponds simply to a 
Schubert cycle, which is irreducible of codimension $(r-r_1)(d_1-r_1)$. If 
we base change $\LG$ over $G_1$ to
$G^{r_1}_1$, we get a proper map with irreducible equidimensional fibers,
mapping surjectively to an irreducible base, so in fact $\LG$ becomes
irreducible, and has dimension precisely $r(d-r)$. Since this dimension 
remains constant
as $r_1$ decreases, and the codimension of $G^{r_1}_1$ increases as $r_1$
decreases, we find we must have exactly one irreducible component of $\LG$
for each choice of $r_1$. 
\end{ex}

\begin{ques}Can we show that $\LG$ is flat over $S$? That it is reduced?
\end{ques}

\begin{ques}Can we describe the components of $\LG$ for $n>2$?
\end{ques}

\bibliographystyle{hamsplain}
\bibliography{hgen}

\newcommand{\noopsort}[1]{} \newcommand{\printfirst}[2]{#1}
  \newcommand{\singleletter}[1]{#1} \newcommand{\switchargs}[2]{#2#1}
\providecommand{\bysame}{\leavevmode\hbox to3em{\hrulefill}\thinspace}
\begin{thebibliography}{10}

\bibitem{b-l-r}
Siegfried Bosch, Werner Lutkebohmert, and Michel Raynaud, \emph{Neron models},
  Springer-Verlag, 1991.

\bibitem{b-l-r2}
Jean-Benoit Bost, Francois Loeser, and Michel Raynaud (eds.), \emph{Courbes
  semi-stables et groupe fondamental en geometrie algebrique}, Birkhauser,
  1998.

\bibitem{b-v}
Winfried Bruns and Udo Vetter, \emph{Determinantal rings}, Lecture Notes in
  Mathematics, no. 1327, Springer-Verlag, 1988.

\bibitem{ci1}
Ciro Ciliberto, \emph{Geometric aspects of polynomial interpolation in more
  variables and of {W}aring's problem}, European Congress of Mathematics, Vol.
  I, Birkhauser, 2001, pp.~289--316.

\bibitem{dj-o1}
A.~J. de~Jong and F.~Oort, \emph{On extending families of curves}, J. Algebraic
  Geometry \textbf{6} (1997), 545--562.

\bibitem{d-m}
P.~Deligne and D.~Mumford, \emph{The irreducibility of the space of curves of
  given genus}, Institut Des Hautes Etudes Scientifiques Publications
  Mathematiques (1969), no.~36, 75--109.

\bibitem{ei1}
David Eisenbud, \emph{Commutative algebra with a view toward algebraic
  geometry}, Springer-Verlag, 1995.

\bibitem{e-h4}
David Eisenbud and Joe Harris, \emph{Divisors on general curves and cuspidal
  rational curves}, Inventiones Mathematicae \textbf{74} (1983), 371--418.

\bibitem{e-h1}
\bysame, \emph{Limit linear series: Basic theory}, Inventiones Mathematicae
  \textbf{85} (1986), 337--371.

\bibitem{es1}
Eduardo Esteves, \emph{Linear systems and ramification points on reducible
  nodal curves}, Matematica Contemporanea \textbf{14} (1998), 21--35,
  \mbox{arXiv:math.AG/9808069}.

\bibitem{gi3}
Alessandro Gimigliano, \emph{Our thin knowledge of fat points}, The Curves
  Seminar at Queen's, Vol. VI, Queen's University, 1989, Exp. No. B.

\bibitem{h-m}
Joe Harris and Ian Morrison, \emph{Moduli of curves}, Springer-Verlag, 1998.

\bibitem{te1}
Montserrat~Teixidor i~Bigas, \emph{Brill-{N}oether theory for stable vector
  bundles}, Duke Mathematical Journal \textbf{62} (1991), no.~2, 385--400.

\bibitem{os7}
B.~Osserman, \emph{Rational functions with given ramification in characteristic
  $p$}, \mbox{arXiv:math.AG/0407445}.

\bibitem{ra1}
Michel Raynaud, \emph{Anneaux locaux henseliens}, Lecture Notes in Mathematics,
  no. 169, Springer-Verlag, 1970.

\bibitem{wi1}
Gayn~B. Winters, \emph{On the existence of certain families of curves},
  American Journal Of Mathematics \textbf{96} (1974), no.~2, 215--228.

\bibitem{ega2}
A.~Grothendieck with J.~Dieudonn\'e, \emph{{\'E}l\'ements de g\'eom\'etrie
  alg\'ebrique: {II.} \'{E}tude globale \'el\'ementaire de quelques classes de
  morphismes}, vol.~8, Publications math\'ematiques de l'I.H.\'E.S., no.~2,
  Institut des Hautes \'Etudes Scientifiques, 1961.

\bibitem{ega32}
\bysame, \emph{{\'E}l\'ements de g\'eom\'etrie alg\'ebrique: {III.} \'{E}tude
  cohomologique des faisceaux coh\'erents, seconde partie}, vol.~17,
  Publications math\'ematiques de l'I.H.\'E.S., no.~2, Institut des Hautes
  \'Etudes Scientifiques, 1963.

\bibitem{ega43}
\bysame, \emph{{\'E}l\'ements de g\'eom\'etrie alg\'ebrique: {IV.} \'{E}tude
  locale des sch\'emas et des morphismes de sch\'emas, troisi\'eme partie},
  vol.~28, Publications math\'ematiques de l'I.H.\'E.S., no.~2, Institut des
  Hautes \'Etudes Scientifiques, 1966.

\bibitem{ega44}
\bysame, \emph{{\'E}l\'ements de g\'eom\'etrie alg\'ebrique: {IV.} \'{E}tude
  locale des sch\'emas et des morphismes de sch\'emas, quatri\'eme partie},
  vol.~32, Publications math\'ematiques de l'I.H.\'E.S., no.~2, Institut des
  Hautes \'Etudes Scientifiques, 1967.

\end{thebibliography}
\end{document}